%% file: 2020_LEPL_SISC.tex
\documentclass[onefignum,onetabnum]{siamart190516}

\input{ex_shared}
\usepackage[caption = false]{subfig}
\usepackage{graphicx}
\usepackage{array}
\newcolumntype{C}[1]{>{\centering\arraybackslash$}p{#1}<{$}}
\ifpdf
\hypersetup{
  pdftitle={Alternating Energy Minimization Methods for Multi-term Matrix Equations},
  pdfauthor={K. Lee, H. C. Elman, C. E. Powell, and D. Lee}
}
\fi

\definecolor{alert}{rgb}{1, 0.12, .75}

\newcommand*{\tran}{^{\mathsf{T}}}
\newcommand*{\fro}{_\text{F}}
\newcommand*{\tolbasis}{\tau_{\text{basis}}}
\newcommand*{\tolblock}{\tau_{\text{coupled}}}

\newcommand{\stagename}{{stage-$p$}}
\newcommand{\stagenameCap}{{Stage-$p$}}
\newcommand{\baselinename}{{S-rank-$1$}}
\newcommand{\pgdaemname}{{PGD-updated}}
\newcommand{\pgdupdatename}{{PGD-update}}
\newcommand{\GSname}{{PGD/Gauss--Seidel-updated}}
\newcommand{\GSnameAbbrv}{{PGD/GS}}
\newcommand{\GSaemnameAbbrv}{{PGD/GS-updated}}
\newcommand{\Rstagename}{{reduced \stagename}}
\newcommand{\RstagenameCap}{{Reduced \stagename}}
\newcommand{\RstagenameAbbrv}{{R-\stagename}}

\newcommand{\ApproxSol}{U_p}

\begin{document}

\maketitle

% REQUIRED
\begin{abstract}
We develop computational methods for approximating the solution of a linear multi-term matrix equation in low rank. We follow an alternating minimization framework, where the solution is represented as a product of two matrices, and approximations to  each matrix are sought by solving certain minimization problems repeatedly. The solution methods we present are based on a rank-adaptive variant of alternating energy minimization methods that builds an approximation iteratively by successively computing a rank-one solution component at each step. We also develop efficient procedures to improve the accuracy of  the low-rank approximate solutions computed using these successive rank-one update techniques. We explore the use of the methods with linear multi-term matrix equations that arise from stochastic Galerkin finite element discretizations of parameterized linear elliptic PDEs, and demonstrate their effectiveness with numerical studies.
\end{abstract}

% REQUIRED

% REQUIRED

\section{Introduction}\label{sec:intro}
We are interested in computing a low-rank approximate solution of a
Kronecker-product structured linear system $Au = b$,
\begin{equation}\label{eq:kron_sys} 
\left(\sum_{i=0}^{m} G_i \otimes K_i \right) u = \sum_{i=0}^{r} g_i \otimes f_i, 
\end{equation} 
where $A =\sum_{i=0}^{m}   G_i \otimes K_i$ is symmetric positive
definite, $\otimes$ is the Kronecker product, $\{K_i\}_{i=0}^{m}\in
\mathbb{R}^{n_1 \times n_1}$, $\{G_i\}_{i=0}^{m} \in \mathbb{R}^{n_2 \times
n_2}$, $\{f_i\}_{i=0}^{r} \in \mathbb{R}^{n_1}$, and $\{g_i\}_{i=0}^{r} \in
\mathbb{R}^{n_2}$. Systems with
such structure arise in the discretization of linear elliptic PDEs in high
dimensions
\cite{ballani2013projection,kressner2015truncated,kressner2010krylov,kressner2011low}
and stochastic Galerkin finite element discretization of parameterized
linear elliptic PDEs
\cite{ghanem2003stochastic,lee2018stochastic,lord2014introduction,
xiu2010numerical}. The solution vector $u \in \mathbb{R}^{n_1 n_2}$ consists of
$n_2$ subvectors of dimension $n_1$, i.e., $u = [u_1\tran, \ldots,
u_{n_2}\tran]\tran$, where $\{u_i\}_{i=1}^{n_2} \in \mathbb{R}^{n_1}$. It also
has an alternative representation in matrix format, $U = [u_1, \ldots, u_{n_2}]
\in \mathbb{R}^{n_1 \times n_2}$, for which the system equivalent to
(\ref{eq:kron_sys}) is the linear multi-term matrix equation \cite{powell2017efficient}
\begin{equation}\label{eq:mat_sys}
\sum_{i=0}^{m} K_i U G_i\tran = B,
\end{equation}
where $B = \sum_{i=0}^{r} f_i g_i\tran \in \mathbb{R}^{n_1 \times n_2}$ and
it is assumed that $m, r \ll n_1, n_2$.  The system matrices $K_i$ and $G_i$ obtained from discretization methods are typically sparse and, thus, for
moderately large system matrices, \textit{Krylov subspace methods}
\cite{pellissetti2000iterative,powell2009block} and \textit{multigrid methods}
\cite{corveleyn2013iterative, elman2007solving,le2003multigrid} have been
natural choices to solve such systems. 

The dimensions of the system matrices grow rapidly, however, if a solution is
sought on a refined grid or (in the case of stochastic Galerkin methods)
if the so-called parametric space is high-dimensional. For large $n_1$ and
$n_2$, direct applications of standard iterative methods may be
computationally prohibitive and storing or explicitly forming the matrix $U$
may be prohibitive in terms of memory. Instead of computing an exact
solution of
\eqref{eq:mat_sys}, we are interested in  inexpensive computation of an
approximate solution of low rank.  To achieve this goal, we begin by
introducing a factored representation of $U \in \mathbb{R}^{n_1 \times n_2}$,
\begin{equation*}
U = VW\tran,
\end{equation*}
where, if $U$ is of full rank $m:=\min(n_1,n_2)$, $V  \in \mathbb{R}^{n_1 \times m}$ and $W  \in \mathbb{R}^{n_2 \times m}$. Our
aim is to find a low-rank approximation to this factored matrix of the form
\begin{equation}\label{eq:apprx_rep}
\ApproxSol = V_{p}W_{p}\tran \in \mathbb{R}^{n_1 \times n_2},
\end{equation}
where $V_p =[v_1,\ldots,v_p] \in \mathbb{R}^{n_1 \times p}$ and $W_p
=[w_1,\ldots,w_p] \in \mathbb{R}^{n_2 \times p}$ and $p \ll m$, and we want to derive
solution algorithms for computing $\ApproxSol$ that operate only on the factors $V_p$
and $W_p$ without explicitly forming $\ApproxSol$.

One such solution algorithm has been developed for matrix completion/sensing \cite{hardt2014understanding,jain2013low}, 
which, at the $p$th iteration,
computes $V_p$ and $W_p$ by alternately solving certain minimization problems.
Although the algorithm computes highly accurate approximations, it can become very expensive as $p$ increases.
Another approach is to use successive rank-one approximations and successively
compute pairs of vectors $\{(v_i$,$w_i)\}_{i=1}^{p}$ to build the factors
$V_p$ and $W_p$ of (\ref{eq:apprx_rep}) until a  stopping criterion is
satisfied. The $p$th iteration starts with $V_{p-1}$ and $W_{p-1}$ and
constructs $v_p$ and $w_p$ as the solutions of certain minimization problems.
This approach for solving parameterized PDEs is one component of a methodology
known as Proper Generalized Decomposition (PGD)
\cite{nouy2007generalized, nouy2010proper,tamellini2014model}. As observed in
those works, using only successive rank-one approximations is
less expensive but may not be practical because it typically results
in approximate solutions with an unnecessarily large value of $p$ for
satisfying a certain error tolerance. 

Our goal in this study is to develop solution algorithms that preserve
only the good properties of the above two types of solution strategies, i.e., algorithms that compute
an accurate solution in a computationally efficient way. In developing
such algorithms, we take our cue from PGD methods, in which, to improve accuracy,
the successive rank-one constructions are supplemented with an {\em updating
procedure} that is performed intermittently during the iteration. Inspired by
this approach, we propose a solution algorithm that adaptively computes 
approximate solutions in an inexpensive way via the successive rank-one
approximation method. This is then supplemented by an enhancement
procedure, which effectively improves the accuracy of the resulting approximate
solutions. We propose two novel enhancement procedures developed by modifying some ideas
used for matrix completion problems \cite{jain2013low}. 

Some other rank-adaptive approaches for approximating solutions of
parameterized or high-dimensional PDEs in low-rank format are as follows.  A
method in \cite{dolgov2014alternating} uses alternating energy minimization
techniques in combination with tensor-train decompositions
\cite{oseledets2011tensor}.  One can incrementally compute rank-one solution
pairs by solving a residual minimization problem, an approach known as
alternating least-squares (ALS) methods, which has been used to compute
low-rank approximate solutions of parameterized PDEs in
\cite{doostan2009least,doostan2013non}, and to solve matrix recovery problems,
matrix sensing and completion problems,
\cite{haldar2009rank,hardt2014understanding,jain2013low,recht2010guaranteed}.
In \cite{powell2017efficient}, an adaptive iterative procedure to solve the matrix equation \eqref{eq:mat_sys} is given, which incrementally
computes a set of orthonormal basis vectors for use in representing the spatial
part of the solution, $V_p$.   See  \cite{simoncini2016computational} for an
overview of other computational approaches for solving linear matrix
equations.

An outline of the paper is as follows.  In Section \ref{sec:aem}, we introduce
and derive alternating energy minimization (AEM) methods using the well-known general
projection framework and discuss a collection of methods developed for
constructing low-rank 
approximate solutions of the form (\ref{eq:apprx_rep}).  In  Section
\ref{sec:update}, we discuss enhancement procedures and derive two new
approaches for performing such updates.  In Section \ref{sec:numexp}, we
measure the effectiveness and the efficiency of the variants of the methods
with numerical experiments.  Finally, in Section \ref{sec:conclusions}, we draw
some conclusions.

\section{Alternating energy minimization (AEM) methods}\label{sec:aem} 
In this section, we derive AEM methods for solving the matrix equation \eqref{eq:mat_sys} from the optimal
projection framework, and review two variants of such methods.  We first introduce some notation. Capital and small
letters are used to denote matrices and vectors, respectively. As a special case, a zero-column matrix is indicated by using a subscript 0, e.g., $X_0\in\mathbb{R}^{n_1 \times 0}$. An inner product between two matrices $X, Y \in \mathbb{R}^{n_1 \times n_2}$ is
defined as $\langle X, Y \rangle \equiv \text{tr}(X\tran Y) = \text{tr}(XY\tran) = \sum_{i,j}
X_{ij}Y_{ij}$, where tr is the trace operator, and $\text{tr}(X) = \sum_{i=1}^{n}
x_{ii}$ if $X\in\mathbb{R}^{n\times n}$. The norm induced by $\langle \cdot, \cdot \rangle$ is the Frobenius norm $\| X \|\fro =
\sqrt{\langle X, X \rangle}$. For shorthand notation, we introduce a linear
operator $\mathcal A( X ) = \sum_{i=0}^{m} K_i X G_i\tran$ for $X \in
\mathbb{R}^{n_1 \times n_2}$. Using this, we can define the weighted inner product
$\langle X, Y\rangle_A = \langle \mathcal A(X), Y\rangle = \langle X, \mathcal
A(Y)\rangle$ and the induced $A$-norm $\| \cdot \|_A$. Finally, vec denotes a vectorization operator, $\text{vec}(X)=x$, where $X =
[x_1,\ldots,x_{n_2}] \in \mathbb{R}^{n_1 \times n_2}$ and
$x=[x_1\tran,\ldots,x_{n_2}\tran]\tran \in \mathbb{R}^{n_1n_2}$, where $x_i \in
\mathbb{R}^{n_1}$,  for $i=1,\ldots,n_2$.

\subsection{General projection framework}\label{sec:prj}
For the computation of $V_p$ and $W_p$ in \eqref{eq:apprx_rep}, we rely on the classical 
theory of orthogonal (Galerkin) projection methods \cite[Proposition
5.2]{saad2003iterative}. Let $\mathcal K \subset \mathbb{R}^{n_1  \times n_2}$ be a \textit{search space} in which an approximate
solution $\ApproxSol \in \mathbb{R}^{n_1 \times n_2}$ is sought, and let
$\mathcal L$ be a \textit{constraint space} onto which the residual $B
-\mathcal A (\ApproxSol)$ is projected. Following \cite[Proposition
5.2]{saad2003iterative}, if the system matrix $A$ is symmetric
positive definite  and $\mathcal L = \mathcal K$, then a matrix $\ApproxSol^\ast$ is
the result of an orthogonal projection onto $\mathcal L$ if and only if
it minimizes the $A$-norm of the error over $\mathcal K$, i.e., 
\begin{equation*}
\ApproxSol^\ast = \underset{\ApproxSol \in \mathcal K}{\arg\min} \:J_A(\ApproxSol),
\end{equation*}
where the objective function is
\begin{equation}\label{eq:Aspd}
J_A(\ApproxSol) = \frac{1}{2} \| U -   \ApproxSol\|_A^2.
\end{equation}
Because we seek a factored representation of $\ApproxSol$, we slightly modify \eqref{eq:Aspd} to give
\begin{equation}\label{eq:opt_obj}
J_A( V_p,W_p) =\frac{1}{2} \| U - V_pW_p\tran \|_A^2, 
\end{equation}
and obtain a new minimization problem 
\begin{equation}\label{eq:min_obj}
\min_{V_p \in \mathbb{R}^{n_1 \times p}, W_p \in \mathbb{R}^{n_2 \times p}} J_A(V_p,W_p).
\end{equation}
Since $J_A$ is quadratic, gradients with respect to $V_p$ and $W_p$ can be easily obtained as
\begin{align}
\nabla_{V_p} J_A &= \left(\mathcal A (V_p W_p\tran) - B \right)W_p = \sum_{i=0}^m (K_i V_p W_p\tran G_i\tran) W_p - BW_p,\label{eq:gradv_Aspd}\\
\nabla_{W_p} J_A &= \left(\mathcal A (V_p W_p\tran) - B \right)\tran V_p = \sum_{i=0}^m  (K_i V_pW_p\tran G_i\tran)\tran V_p  -  B\tran V_p.\label{eq:gradw_Aspd}
\end{align}
Employing the first-order optimality condition on \eqref{eq:gradv_Aspd}--\eqref{eq:gradw_Aspd} 
 (i.e., setting \eqref{eq:gradv_Aspd} and (the transpose of) \eqref{eq:gradw_Aspd}
to be zero) results in the set of equations
\begin{align}
\sum_{i=0}^m (K_i V_p W_p\tran G_i\tran) W_p &= BW_p \in \mathbb{R}^{n_1 \times  p},\label{eq:vp_Aspd}\\
\sum_{i=0}^m V_p\tran (K_i V_p W_p\tran G_i \tran)  &= V_p\tran B \in \mathbb{R}^{ p \times n_2 } .\label{eq:wq_Aspd}
\end{align}
These equations can be interpreted as projections of the residual $B - \mathcal A (V_p W_p\tran) $
onto the spaces spanned by the columns of $W_p$ and $V_p$, respectively. 

Given \eqref{eq:vp_Aspd}--\eqref{eq:wq_Aspd}, a
widely used strategy for solving the minimization problem
\eqref{eq:min_obj} is to compute each component of the solution pair
$(V_p,W_p)$ alternately
\cite{dolgov2014alternating,doostan2009least,doostan2013non,haldar2009rank,hardt2014understanding,jain2013low}.
That is, one can fix $W_p$ and solve the system of equations of order $n_1  p$
in \eqref{eq:vp_Aspd}  for $V_p$, and then one can fix $V_p$ and solve the
system of equations of order $n_2  p$ in \eqref{eq:wq_Aspd}  for  $W_p$.  However, in this approach, suitable choices of $p$ for satisfying a fixed
error tolerance are typically not known \textit{a priori}. Thus, adaptive schemes
that incrementally compute solution pairs ($v_i$,$w_i$) have been introduced
\cite{jain2013low,nouy2007generalized, nouy2010proper,tamellini2014model}. All
of these schemes are based on alternately solving two systems of equations for
two types of variables in an effort to minimize a certain error measure. In
this study, we employ alternating methods for minimizing the
energy norm of the error \eqref{eq:min_obj} and,
thus, we refer to approaches of this type as \textit{alternating energy
minimization} (AEM) methods. In the following sections, we present two 
adaptive variants of AEM methods: a
\stagenameCap\ AEM method 
and a successive rank-one AEM method. 

\subsection{\stagenameCap\ AEM method}\label{sec:stage}
An alternating minimization method that entails solving a sequence of 
least-squares problems whose dimensions increase with $p$ was developed in \cite{jain2013low} 
for solving matrix-recovery problems \cite{haldar2009rank,hardt2014understanding,jain2013low}.
We adapt this approach to the energy
minimization problem \eqref{eq:min_obj} and refer to it as the \stagenameCap\  AEM 
method.  
It is an iterative method that runs until an approximate solution
satisfies a  stopping criterion (e.g., the relative residual $\|B - \mathcal
A(V_pW_p\tran)\|\fro \leq \epsilon \| B\|\fro$ with a user-specified
stopping tolerance $\epsilon$). 
At the $p$th iteration, 
called a ``stage'' in \cite{jain2013low}, this method seeks 
$p$-column factors $V_p$ and $W_p$
determining an approximate solution by initializing $W_p^{(0)}$ and solving the following
systems of equations in sequence:
\begin{alignat}{2}
\sum_{i=0}^m (K_i) V_p^{(k)} (W_p^{(k-1)}{}\tran G_i W_p^{(k-1)})\tran &= BW_p^{(k-1)}, \label{eq:vp_stage}\\ 
\sum_{i=0}^m (V_p^{(k)}{}\tran K_i V_p^{(k)}) W_p^{(k)}{}\tran  (G_i \tran)  &= V_p^{(k)}{}\tran B, \label{eq:wp_stage}
\end{alignat}
for $k=1,\ldots,k_{\max}$, where the superscript indicates the number of
alternations between the two systems of equations
\eqref{eq:vp_stage}--\eqref{eq:wp_stage}. Note that the method can also begin by
initializing $V_p^{(0)}$ and alternating between
\eqref{eq:wp_stage} and \eqref{eq:vp_stage}. Algorithm \ref{alg:stage_aem}
summarizes the entire procedure.  For the initialization of $W_p^{(0)}$ (line
\ref{alg:svp}), one step of the singular value projection method
\cite{jain2010guaranteed} is performed with the exact settings from
\cite[Algorithm 3]{jain2013low}.  The \textsc{CheckConvergence} procedure 
(line \ref{alg:check_conv_stagep})
is detailed in Section \ref{sec:update}.

\begin{algorithm}[!h]
\caption{\stagenameCap\  AEM method} % (stagewise-AEM)}
\hspace*{\algorithmicindent} \textbf{INPUT:} 
%$B$: the right-hand side of the matrix equation,\\
$p_{\max}$: the maximum number of solution pairs, \\
\hspace*{21mm} $k_{\max}$: the maximum number of alternations in each stage, \\
\hspace*{21mm} $\epsilon$: a parameter for checking convergence, 
\begin{algorithmic}[1]
\Function{StagepAEM}{$p_{\max}, k_{\max}, \epsilon$}
  \For{ $p=1,\ldots,p_{\max}$}
    \State $[V_p^{(0)}, W_p^{(0)}] =$ first $p$ singular vectors of $V_{p-1} W_{p-1}\tran - \frac{3}{4}(\mathcal A(V_{p-1} W_{p-1}\tran) \!-\! B)$    \label{alg:svp}  
    \For{ $k=1,\ldots,k_{\max}$}
      \State $V_p^{(k)}  \leftarrow$ solve \eqref{eq:vp_stage} \label{ll:Vp} % \underset{{V_p \in \mathbb{R}^{n_1\times p}} }{\arg \min}\: J_A(V_p, W_p^{(k-1)})$
      \State $W_p^{(k)}  \leftarrow$ solve \eqref{eq:wp_stage}  \label{ll:Wp}% \underset{{W_p \in \mathbb{R}^{n_2 \times p}} }{\arg \min}\: J_A(V_p^{(k)},W_p)$
    \EndFor
    \State $V_p \leftarrow V_p^{(k)}$ and $W_p \leftarrow W_p^{(k)}$
    %\If{$\|B - \mathcal A(V_pW_p\tran) \|\fro < \epsilon \|B\|\fro$}
    % \State break
    %\EndIf
    \State $V_p, W_p \leftarrow $ \Call{CheckConvergence}{$V_p, W_p,\epsilon$} \label{alg:check_conv_stagep}
  \EndFor
\EndFunction
\end{algorithmic}
\label{alg:stage_aem}
\end{algorithm}

Systems of equations for ``vectorized'' versions of the matrix factors
$V_p$ and $W_p$ can be derived\footnote{The left-hand sides of
\eqref{eq:stage_vp}--\eqref{eq:stage_wp} are derived using $\text{vec}(KUG\tran) =
(G\otimes K) \text{vec}(U)$. Note that \eqref{eq:stage_wp} is derived
by first transposing \eqref{eq:wp_stage} and then vectorizing the resulting
equation. In the sequel, vectorized versions of equations for the factor
$W_p$ are derived by first taking the transpose.} from \eqref{eq:vp_stage} and
\eqref{eq:wp_stage} as follows
\begin{align}
\sum_{i=0}^{m} [ (W_p^{(k-1)}{}\tran G_i W_p^{(k-1)})\otimes K_i ]\,\text{vec}(V_p^{(k)}) &= \text{vec}(B W_p^{(k-1)}),\label{eq:stage_vp}\\
\sum_{i=0}^{m} [ (V_p^{(k)}{}\tran K_i V_p^{(k)}) \otimes G_i ]\,\text{vec}(W_p^{(k)}) & = \text{vec}( B\tran V_p^{(k)}{}).\label{eq:stage_wp}
\end{align}
Thus, solving \eqref{eq:vp_stage} and \eqref{eq:wp_stage} is equivalent to
solving coupled linear systems with coefficient matrices of dimensions $n_1 p \times n_1 p$ and $n_2 p
\times n_2 p$, respectively, which are smaller than that of the original system
\eqref{eq:mat_sys} when $p$ is small. However, the reduced matrix factors (of size $p \times p$) are dense, even if the original ones are sparse, and so as $p$ increases,
the computational costs for solving \eqref{eq:vp_stage}--\eqref{eq:wp_stage} increase and the 
\stagenameCap\ AEM %stagewise-AEM
method may be impractical for large-scale problems.

\subsection{Successive rank-one AEM method}\label{sec:succ_rk1}
We now describe a successive rank-one (\baselinename) approximation method which, at each
iteration, adds a rank-one correction to the current iterate. This is
a basic component of PGD methods \cite{nouy2007generalized,
nouy2010proper,tamellini2014model} for solving parameterized PDEs. The method
only requires solutions of linear systems with coefficient matrices of size 
$n_1 \times n_1$ and $n_2 \times n_2$ rather than coupled systems 
like those in the \stagenameCap\ AEM method that grow in size with the step
counter $p$. % \CP{as in the stagewise-AEM method.}

Assume that $p-1$ pairs of solutions are computed, giving $V_{p-1}$ and
$W_{p-1}$.  The next step is to compute a new solution pair $(v_p,w_p)$ by choosing
the objective function
\begin{equation*}
J_A( v_p, w_p) =\frac{1}{2} \| U - V_{p-1}W_{p-1}\tran - v_p w_p\tran \|_A^2,
\end{equation*}
and solving the following minimization problem
\begin{equation*}
\min_{v_p \in \mathbb{R}^{n_1}, w_p \in \mathbb{R}^{n_2}} J_A(v_p, w_p).
\end{equation*}
The gradients of $J_A$ with respect to $v_p$ and $w_p$ are
\begin{align}
\nabla_{v_p} J_A &= \left( \mathcal A(v_{p} w_{p}\tran) + \mathcal A(V_{p-1} W_{p-1}\tran) - B\right)w_p,\label{eq:gradv_Aspd_rk1}\\  
\nabla_{w_p} J_A &= \left( \mathcal A(v_{p} w_{p}\tran) + \mathcal A(V_{p-1} W_{p-1}\tran) - B\right)\tran v_p.\label{eq:gradw_Aspd_rk1}
\end{align}

Employing the first-order optimality conditions
(setting \eqref{eq:gradv_Aspd_rk1} and (the transpose of) \eqref{eq:gradw_Aspd_rk1} to zero)
results in systems of equations for
which, in a succession of steps $k=1,\ldots,k_{\max}$, $v_p$ is updated using
fixed $w_p$ and then $w_p$ is updated using fixed $v_p$:
\begin{align}
\sum_{i=0}^m (K_i) v_p^{(k)} (w_p^{(k-1)}{}\tran G_i w_p^{(k-1)})\tran  &= Bw_p^{(k-1)} - \mathcal A( V_{p-1}W_{p-1}\tran) w_p^{(k-1)}, \label{eq:vp_rk1}\\
 \sum_{i=0}^m (v_p^{(k)}{}\tran K_i  v_p^{(k)}) w_p^{(k)}{}\tran (G_i\tran)  &= v_p^{(k)}{}\tran B -  v_p^{(k)}{}\tran \mathcal A( V_{p-1} W_{p-1}\tran). \label{eq:wp_rk1}
\end{align}
Algorithm \ref{alg:base_aem} summarizes this procedure, which randomly
initializes $w_p^{(0)}$ and then alternately solves
\eqref{eq:vp_rk1}--\eqref{eq:wp_rk1}.
Like the \stagenameCap\ AEM method, the
algorithm can start with
either $w_p^{(0)}$ or $v_p^{(0)}$.
\begin{algorithm}[!t]
\caption{Successive rank-one AEM method} % (baseline-AEM)}
\hspace*{\algorithmicindent} \textbf{INPUT:} $p_{\max}$, $k_{\max}$, and $\epsilon$
\begin{algorithmic}[1]
\Function{SrankoneAEM}{$p_{\max}$, $k_{\max}$, $\epsilon$}
\For{ $p=1,\ldots,p_{\max}$}
  \State Set a random initial guess for $w_p^{(0)}$.
  \For{ $k=1,\ldots,k_{\max}$}
  \State $v_p^{(k)}  \leftarrow$ solve \eqref{eq:vp_rk1} 
  \State $w_p^{(k)}  \leftarrow$ solve \eqref{eq:wp_rk1} 
  \EndFor
  \State $v_p \leftarrow v_p^{(k)}$ and $w_p \leftarrow w_p^{(k)}$
  \State Add to solution matrices, $V_p \leftarrow [V_{p-1}, v_p]$, $W_p \leftarrow [W_{p-1}, w_p]$
  \State $V_p, W_p \leftarrow $ \Call{CheckConvergence}{$V_p, W_p,\epsilon$}
\EndFor
\EndFunction
\end{algorithmic}
\label{alg:base_aem}
\end{algorithm}

\subsection{Algebraic interpretation of the methods}\label{sec:interpretation}
%It can be seen that
%both methods, stagewise-AEM (Algorithm \ref{alg:stage_aem})
%and baseline-AEM (Algorithm \ref{alg:base_aem}), 
Algorithms \ref{alg:stage_aem} and \ref{alg:base_aem} both
entail an ``outer iteration''
with counter $p$ and an ``inner iteration'' with counter $k$, 
%to achieve the same goal: find a minimizer of 
and both are designed to minimize
the objective function \eqref{eq:opt_obj}. 
%It will be 
It is instructive to see the difference between the two methods in vectorized
format. To this end, let 
 \begin{align*}
{\mathcal A_w}(w_i,w_j) = \sum_{{l=0}}^m K_{l} (w_j\tran G_{l}\tran w_i) \in \mathbb{R}^{n_1 \times n_1},\quad%\label{eq:v_op}\\ 
{\mathcal A_v}(v_i,v_j) = \sum_{{l=0}}^m G_{l} (v_j\tran K_{l}\tran v_i) \in \mathbb{R}^{n_2 \times n_2},%\label{eq:w_op}
\end{align*}
and let us assume $p=2$ for simplifying the presentation. 

%Both methods seek solution pairs $(V_3, W_3)$ satisfying the systems of
Both methods seek solution pairs $(V_2, W_2)$ satisfying the systems of
equations \eqref{eq:vp_Aspd}--\eqref{eq:wq_Aspd}, which can be written in a
vectorized form:
\begin{align}
\left[
\begin{array}{*4{C{1.8cm}}}
A_w(w_1,w_1) & A_w(w_1,w_2)\\ %& A_w(w_1,w_3)\\
A_w(w_2,w_1) & A_w(w_2,w_2)\\ %& A_w(w_2,w_3)\\
%A_w(w_3,w_1) & A_w(w_3,w_2) & A_w(w_3,w_3)\\
\end{array}
\right]
\left[
\begin{array}{*4{C{.4cm}}}
v_1 \\
v_2 %\\
%v_3
\end{array}
\right]
&= 
\left[
\begin{array}{*4{C{.8cm}}}
B w_1\\
B w_2%\\
%B w_3
\end{array}
\right], \label{eq:vp_Aspd_vec}
\\
\left[
\begin{array}{*4{C{1.8cm}}}
A_v(v_1,v_1) & A_v(v_1,v_2) \\%& A_v(v_1,v_3)\\
A_v(v_2,v_1) & A_v(v_2,v_2) \\%& A_v(v_2,v_3)\\
%A_v(v_3,v_1) & A_v(v_3,v_2) & A_v(v_3,v_3)\\
\end{array}
\right]
\left[
\begin{array}{*4{C{.4cm}}}
w_1 \\
w_2 %\\
%w_3
\end{array}
\right]
&= 
\left[
\begin{array}{*4{C{.8cm}}}
B\tran v_1\\
B\tran v_2%\\
%B\tran v_3
\end{array}
\right].\label{eq:wp_Aspd_vec}
\end{align}
%In the third outer iteration, the  \stagename\ AEM method alternately solves
In the second outer iteration, the  \stagenameCap\ AEM method alternately solves
fully coupled linear systems \eqref{eq:vp_stage}--\eqref{eq:wp_stage}
%specified by $W_3^{(k-1)}$ and $V_3^{(k)}$, which can be written in 
specified by $W_2^{(k-1)}$ and $V_2^{(k)}$, which can be written in 
vectorized form as in \eqref{eq:vp_Aspd_vec}--\eqref{eq:wp_Aspd_vec}:
\begin{align}
\left[
\begin{array}{*4{C{2.8cm}}}
\!\!\!A_w(w_1^{(k-1)},w_1^{(k-1)})\!\! & \!\!\!A_w(w_1^{(k-1)},w_2^{(k-1)})\!\!\\ %& \!\!\!A_w(w_1^{(k-1)},w_3^{(k-1)})\!\!\!\\
\!\!\!A_w(w_2^{(k-1)},w_1^{(k-1)})\!\! &\!\!\! A_w(w_2^{(k-1)},w_2^{(k-1)})\!\!\\ %&\!\!\!A_w(w_2^{(k-1)},w_3^{(k-1)})\!\!\!\\
%\!\!\!A_w(w_3^{(k-1)},w_1^{(k-1)}) \!\!& \!\!\!A_w(w_3^{(k-1)},w_2^{(k-1)}) \!\!& \!\!\!A_w(w_3^{(k-1)},w_3^{(k-1)})\!\!\!\\
\end{array}
\right]\!\!\!
\left[
\begin{array}{*4{C{.3cm}}}
\!\!\!v_1^{(k)} \!\!\!\\
\!\!\!v_2^{(k)} \!\!\!%\\
%\!\!\!v_3^{(k)} \!\!\!
\end{array}
\right]
&\!\!= \!\!
\left[
\begin{array}{*4{C{0.9cm}}}
\!\!\!\!B w_1^{(k-1)}\!\!\!\!\\
\!\!\!\!B w_2^{(k-1)}\!\!\!\!%\\
%\!\!\!\!B w_3^{(k-1)}\!\!\!\!
\end{array}
\right],\!\! \nonumber
\\
\left[
%\begin{array}{*4{C{2.35cm}}}
\begin{array}{*4{C{2.8cm}}}
\!\!A_v(v_1^{(k)},v_1^{(k)}) \!\!&\!\! A_v(v_1^{(k)},v_2^{(k)}) \!\!\\%&\!\! A_v(v_1^{(k)},v_3^{(k)})\!\!\\
\!\!A_v(v_2^{(k)},v_1^{(k)}) \!\!&\!\! A_v(v_2^{(k)},v_2^{(k)}) \!\!\\%&\!\! A_v(v_2^{(k)},v_3^{(k)})\!\!\\
%\!\!A_v(v_3^{(k)},v_1^{(k)}) \!\!&\!\! A_v(v_3^{(k)},v_2^{(k)}) \!\!&\!\! A_v(v_3^{(k)},v_3^{(k)})\!\!\\
\end{array}
\right]\!\!\!
\left[
\begin{array}{*4{C{.3cm}}}
\!\!\!\!w_1^{(k)}\!\!\!\! \\
\!\!\!\!w_2^{(k)}\!\!\!\! %\\
%\!\!\!\!w_3^{(k)}\!\!\!\!
\end{array}
\right]
&\!\!= \!\!
\left[
\begin{array}{*4{C{0.9cm}}}
\!\!B\tran v_1^{(k)}\!\!\\
\!\!B\tran v_2^{(k)}\!\!%\\
%\!\!B\tran v_3^{(k)}\!\!
\end{array}\!\!
\right]. \label{eq:wp_stagevec}
\end{align}

In contrast, the \baselinename\ method seeks approximate solutions of
\eqref{eq:vp_Aspd_vec}--\eqref{eq:wp_Aspd_vec} by solving systems of equations
associated with only the diagonal blocks. %of \eqref{eq:vp_Aspd_vec}--\eqref{eq:wp_Aspd_vec}. 
%In the first outer iteration, the \baselinename\ method alternates between the following two equations to find
In the first iteration, the method alternates between the following equations to find
$v_1$ and $w_1$:
\begin{align*}
\left[
\begin{array}{*4{C{2.7cm}}}
\!\!\!A_w(w_1^{(k-1)},w_1^{(k-1)})\!\!\! & %&
\end{array}
\right]\!\!\!
\left[
\begin{array}{*4{C{.35cm}}}
\!\!v_1^{(k)} \!\!\!
\end{array}
\right]
&\!\!= \!\!
\left[
\begin{array}{*4{C{1.1cm}}}
\!\!\!B w_1^{(k-1)}\!\!\!
\end{array}
\right],
\\
\left[
\begin{array}{*4{C{2.7cm}}}
\!\!A_v(v_1^{(k)},v_1^{(k)}) \!\! & %&
\end{array}
\right]\!\!\!
\left[
\begin{array}{*4{C{.35cm}}}
\!\!\!w_1^{(k)}\!\!\! 
\end{array}
\right]
&\!\!= \!\!
\left[
\begin{array}{*4{C{1.1cm}}}
B\tran v_1^{(k)}
\end{array}
\right].
\end{align*}
In the second iteration, the method alternately solves the systems of equations
in the second rows of the following equations:
\begin{align*}
\left[
\begin{array}{*4{C{2.7cm}}}
\!\!\!A_w(w_1,w_1)\!\!\! &  \\%& \\
\!\!\!A_w(w_2^{(k-1)},w_1)\!\!\! &\!\!\! A_w(w_2^{(k-1)},w_2^{(k-1)})\!\!\! %&
\end{array}
\right]\!\!\!
\left[
\begin{array}{*4{C{.35cm}}}
\!\!v_1 \!\!\!\\
\!\!v_2^{(k)} \!\!\!
\end{array}
\right]
&\!\!= \!\!
\left[
\begin{array}{*4{C{1.1cm}}}
\!\!\!B w_1\!\!\!\\
\!\!\!B w_2^{(k-1)}\!\!\!
\end{array}
\right],
\\
\left[
\begin{array}{*4{C{2.7cm}}}
\!\!A_v(v_1,v_1) \!\!& \\%&\\
\!\!A_v(v_2^{(k)},v_1) \!\!&\!\! A_v(v_2^{(k)},v_2^{(k)}) \!\!%&
\end{array}
\right]\!\!\!
\left[
\begin{array}{*4{C{.35cm}}}
\!\!\!w_1\!\!\! \\
\!\!\!w_2^{(k)}\!\!\! 
\end{array}
\right]
&\!\!= \!\!
\left[
\begin{array}{*4{C{1.1cm}}}
B\tran v_1\\
B\tran v_2^{(k)}
\end{array}
\right].
\end{align*}
Because $v_1$ and $w_1$ are fixed, the (2,1)-block matrices are multiplied with
$v_1$ and $w_1$ and the resulting vectors are moved to the right-hand sides.
Then solving the equations associated with the (2,2)-block matrices gives
$v_2^{(k)}$ and $w_2^{(k)}$. 
As illustrated in this example, the
\baselinename\ AEM method approximately solves \eqref{eq:vp_Aspd_vec}--\eqref{eq:wp_Aspd_vec} by taking the
matrices in the lower-triangular blocks to the right-hand sides and solving only the
systems associated with the diagonal blocks, as opposed to solving fully coupled systems as in
the \stagenameCap\ AEM method.

The system matrices that arise in Algorithm \ref{alg:stage_aem} have reduced factors that are dense but small (of size $p \times p$) and their counterpart factors are large but sparse. In Algorithm \ref{alg:base_aem}, the system matrices are sparse and of order $n_1$ and $n_2$ (as the reduced factors are of size 1 $\times$ 1). Thus in both cases, we may use Krylov subspace methods to solve the systems. Then, with the iteration counter $p$, the cost of the \stagenameCap\ AEM method grows quadratically (since the reduced factors are dense), whereas that of the
\baselinename\ AEM method grows linearly with $p$. Thus, using the \stagenameCap\ AEM method can be impractical for large-scale
applications. On the other hand, as the \baselinename\ AEM method employs only
the lower-triangular part of the system matrices, convergence tends to be
slow and the level of accuracy that can be achieved in a
small number of steps is limited. 
To overcome these shortcomings, 
in the next section, we will consider several ways to modify and enhance them to improve accuracy.

\remark{The \stagenameCap\ AEM and \baselinename\ AEM methods can be seen as two extreme versions of AEM methods. The former solves fully coupled systems and the latter sequentially solves systems associated with the diagonal blocks. Although it has not been explored in this study, in an intermediate approach, more than one consecutive pair of solution vectors $(\{v_{p},\ldots,v_{p+\ell}\},\{w_{p},\ldots,w_{p+\ell}\})$, with $\ell \in \mathbb{N}$, can be computed in a coupled manner at each outer iteration.}

\section{Enhancements}\label{sec:update}
We now describe variants of the \baselinename\,AEM method 
that perform extra computations to improve accuracy. 
The general strategy is to compute an enhancement of the approximate solution
at every $n_{\text{update}}$ outer iterations of the \baselinename\ AEM method, as specified in Algorithms
\ref{alg:enhancedAEM}--\ref{alg:check_conv}. 

\begin{algorithm}[!h]
\caption{Enhanced AEM method}
\hspace*{\algorithmicindent} \textbf{INPUT:} $p_{\text{max}}$, $k_{\max}$, $n_{\text{update}}$, and $\epsilon$
\begin{algorithmic}[1]
\Function{EnhancedAEM}{$p_{\max}, k_{\max}, n_{\text{update}}, \epsilon$}
\For{ $p=1,\ldots,p_{\max}$}
\State $v_p, w_p \leftarrow$ \Call{RankOneCorrection}{$V_{p-1},W_{p-1},k_{\max}$}
\State Add to solution matrices, $V_p \leftarrow [V_{p-1}, v_p]$, $W_p \leftarrow [W_{p-1}, w_p]$
\If {$p \mod n_{\text{update}}==0$} \label{line:nupdate}
\State $V_p, W_p \leftarrow$ \Call{Enhancement}{$V_p, W_p$}
\EndIf
\State $V_p, W_p \leftarrow $ \Call{CheckConvergence}{$V_p, W_p,\epsilon$}
\EndFor
\EndFunction
\end{algorithmic}
\label{alg:enhancedAEM}
\end{algorithm}

\begin{algorithm}[!h]
\caption{Rank one correction}
\hspace*{\algorithmicindent} \textbf{INPUT:} $V_{p-1},W_{p-1}$, and $k_{\max}$
\begin{algorithmic}[1]
\Function{RankOneCorrection}{$V_{p-1}$,$W_{p-1}$,$k_{\max}$}
\State Set a random initial guess for $w_p^{(0)}$.
\For{ $k=1,\ldots,k_{\max}$}
\State $v_p^{(k)}  \leftarrow$ solve \eqref{eq:vp_rk1}
\State $w_p^{(k)}  \leftarrow$ solve \eqref{eq:wp_rk1} 
\EndFor
\State $v_p \leftarrow v_p^{(k)}$ and $w_p \leftarrow w_p^{(k)}$
\EndFunction
\end{algorithmic}
\label{alg:rank1_correction}
\end{algorithm}

\begin{algorithm}[!h]
\caption{Checking for convergence}
\hspace*{\algorithmicindent} \textbf{INPUT:} $V_{p},W_{p}$, and $\epsilon$
\begin{algorithmic}[1]
\Function{CheckConvergence}{$V_{p}$,\,$W_{p}$,\,$\epsilon$}
  \If{$\|V_p W_p\tran - V_{p-1}W_{p-1}\tran\|\fro \leq \epsilon \| V_pW_p\tran\|\fro$}
    \State $V_p, W_p \leftarrow$ \Call{Enhancement}{$V_p, W_p$}
    \If{$\|V_p W_p\tran - V_{p-1}W_{p-1}\tran\|\fro \leq \epsilon \| V_pW_p\tran\|\fro$} Stop
    \EndIf
  \EndIf
\EndFunction
\end{algorithmic}
\label{alg:check_conv}
\end{algorithm}

We present three enhancement procedures, 
one taken from the literature and two new ones.
These are 
(i) a procedure adopted from an updating technique developed in 
\cite[Section 2.5]{tamellini2014model}, which defines one variant of PGD methods; 
(ii) a refined version of this approach,  which only solves systems 
associated with the diagonal blocks of the system matrices but incorporates
information (upper-triangular blocks)
in a manner similar to Gauss-Seidel iterations;
and (iii) an adaptive  enhancement of the \stagenameCap\ AEM method that 
decreases costs with negligible impact on accuracy.
In discussing these ideas, we distinguish updated solutions
using the notation, $\overline{v}_i$, $\overline{w}_i$ (for vectors), and
$\overline{V}_p = [\overline{v}_1,\ldots,\overline{v}_p]$, $\overline{W}_p =
[\overline{w}_1,\ldots,\overline{w}_p]$ (for matrices).

Before we detail each method, we first elaborate on the $\textsc{CheckConvergence}$ procedure in Algorithm
\ref{alg:check_conv}. This checks the relative difference
between the current iterate and the previous iterate $\|V_pW_p\tran - V_{p-1}W_{p-1}\tran\|\fro \leq \epsilon \|V_pW_p\tran\|\fro$  in the Frobenius
norm.\footnote{To compute $ \|V_pW_p\tran\|\fro^{2}$, we form $X = (V_p\tran
V_p) \odot (W_p\tran W_p) \in \mathbb{R}^{p\times p}$, where $\odot$ is the Hadamard product, and then sum-up all the elements of X. The product $V_p W_p\tran$ is never explicitly formed.} If this condition is met, we apply the $\textsc{Enhancement}$
procedure and check the convergence with
the same criterion. The purpose of this extra enhancement %\textsc{Enhancement} 
is to help prevent
Algorithm \ref{alg:enhancedAEM} from terminating prematurely (i.e., the
stopping condition can be met when Algorithm \ref{alg:enhancedAEM}
stagnates.).

%\subsection{PGD update / single-sided enhancement methods}\label{sec:sr}

\subsection{\pgdaemname\ AEM}\label{sec:sr}
Suppose the factors $V_p$ and $W_p$ obtained from
\textsc{Rank\-One\-Correction} do not satisfy the first-order optimality conditions
\eqref{eq:vp_Aspd}--\eqref{eq:wq_Aspd}. 
%Following the PGD updating procedure
An enhancement like that of the PGD update
\cite{nouy2007generalized, nouy2010proper,tamellini2014model} 
%this enhancement procedure 
modifies one of these factors (e.g., the one corresponding to
the smaller dimension $n_1$ or $n_2$) by solving the associated
minimization problem for $V_p$ (given $W_p$, when $n_1<n_2$) or for $W_p$
(given $V_p$ when $n_1>n_2$) so that one of the first-order conditions holds.
We outline the procedure for approximating $W_p$;
the procedure for $V_p$ is analogous. %approximating $V_p$ is entirely analogous.  
The basic procedure is to solve the optimization problem 
$\min_{W_p \in \mathbb{R}^{n_2 \times p}} J_A\left( V_p, W_p\right)$
every $n_{\text update}$ steps. In place of $V_p$, an orthonormal matrix 
${\widetilde V}_p$ is used, so that the construction entails solving
\begin{equation}\label{eq:os_obj}
\overline{W}_p = \underset{W_p \in \mathbb{R}^{n_2 \times p}}{\arg\min}J_A\left( \widetilde V_p,  W_p\right),
\end{equation}
where $J_A$ is the quadratic objective function defined in \eqref{eq:opt_obj}.
The gradient of the objective function $J_A$ with respect to  $W_p$ can be
computed as
\begin{align*}
\nabla_{W_p} J_A &= \left(\mathcal A (\widetilde V_p W_p\tran) - B \right)\tran\widetilde V_p = \sum_{i=0}^m  (K_i \widetilde V_p W_p\tran G_i\tran)\tran \widetilde V_p  - B\tran\widetilde V_p.
\end{align*}
Thus, solving the minimization problem \eqref{eq:os_obj} by employing the
first-order optimality condition is equivalent to solving a system of equations %,
%which is analogous to  \eqref{eq:wq_Aspd},
similar in structure to  \eqref{eq:wq_Aspd},
\begin{align}
\sum_{i=0}^m (\widetilde V_p\tran K_i \widetilde V_p) \overline{W}_p\tran  (G_i \tran)  &= \widetilde V_p\tran B \in \mathbb{R}^{ p \times n_2 }.\label{eq:wq_Aspd_os}
\end{align}

Compared to the original system \eqref{eq:mat_sys}, the dimension of this
matrix is reduced via a ``single-sided'' reduction; in
\eqref{eq:wq_Aspd_os}, the reduction is on the side of the
first dimension, i.e., $n_1$ is reduced to $p$. The vectorized form of this system, for
$p=2$, is
\begin{align*}
\left[
\begin{array}{*4{C{2.7cm}}}
\!\!A_v(\tilde v_1,\tilde v_1) \!\!&\!\! A_v(\tilde v_1,\tilde v_2) \!\!\\%&\!\! A_v(\tilde v_1,\tilde v_3)\!\!\\
\!\!A_v(\tilde v_2,\tilde v_1) \!\!&\!\! A_v(\tilde v_2,\tilde v_2) \!\!\\%&\!\! A_v(\tilde v_2,\tilde v_3)\!\!\\
%\!\!A_v(\tilde v_3,\tilde v_1) \!\!&\!\! A_v(\tilde v_3,\tilde v_2) \!\!&\!\! A_v(\tilde v_3,\tilde v_3)\!\!\\
\end{array}
\right]\!\!\!
\left[
\begin{array}{*4{C{.3cm}}}
\!\!\!\overline{w}_1\!\!\! \\
\!\!\!\overline{w}_2\!\!\! %\\
%\!\!\!\overline{w}_3\!\!\!
\end{array}
\right]
&\!\!= \!\!
\left[
\begin{array}{*4{C{.8cm}}}
B\tran \tilde v_1\\
B\tran \tilde v_2%\\
%B\tran \tilde v_3
\end{array}
\right],
\end{align*}
which has structure like that of the second system in \eqref{eq:wp_stagevec} of the
\stagenameCap\, AEM method. 
We summarize this single-sided enhancement method in Algorithm \ref{alg:aem_sr}.

\remark{Another approach for computing a set of orthonormal basis vectors and computing a low-rank solution by solving a reduced system of type \eqref{eq:wq_Aspd_os} is given in \cite{powell2017efficient}. The MultiRB method of \cite{powell2017efficient} incrementally computes a set of orthonormal basis vectors for the spatial part of the solution (i.e., $\widetilde V_p \in \mathbb{R}^{n_1 \times p}$) using \textit{rational Krylov subspace methods} and solves a reduced system for $\overline{W}_p$ and, consequently, $\ApproxSol = \widetilde V_p \overline{W}_p\tran$.}

\begin{algorithm}[!t]
\caption{\pgdupdatename\ enhancement}
\hspace*{\algorithmicindent} \textbf{Input:} $V_p$ and $W_p$
\begin{algorithmic}[1]
\Function{PGDupdate}{$V_p,W_p$}
\If{$n_1 < n_2$}
\State $ \widetilde W_p \leftarrow$ orthonormalize $W_p$.
\State ${\overline{V}_p} \leftarrow $ solve $\sum_{i=0}^m (K_i) \overline{V}_p (\widetilde{W}_p\tran G_i \widetilde{W}_p)\tran  = B\widetilde{W}_p$
\State {$V_p \leftarrow \overline{V}_p$}
\Else
\State $\widetilde V_p \leftarrow $ orthonormalize $V_p$.
\State ${\overline{W}_p} \leftarrow$ solve $\sum_{i=0}^m (\widetilde V_p\tran K_i \widetilde V_p) \overline{W}_p\tran  (G_i \tran)  = \widetilde V_p\tran B$
\State {$W_p \leftarrow \overline{W}_p$}
\EndIf
\EndFunction
\end{algorithmic}
\label{alg:aem_sr}
\end{algorithm}

\subsection{\GSname\ AEM} \label{sec:diag}
The second strategy for enhancement,
like the ``unenhanced'' \baselinename\ AEM method  
(and in contrast to \pgdaemname\ AEM), 
only requires solutions of linear systems with coefficient matrices of 
dimensions $n_1 \times n_1$ and $n_2 \times n_2$, independent of $p$.
As observed in Section
\ref{sec:interpretation}, the \baselinename\ AEM method loosely corresponds to
solving lower block-triangular systems of equations.  We modify these
computations by using more information 
(from the upper triangular part), as soon
as it becomes available.  This leads to a method that resembles
the (block) Gauss--Seidel method for linear systems \cite{golub2012matrix}.
Suppose $\{(v_i, w_i)\}_{i=1}^{p}$ are obtained from $p$ iterations of
Algorithm \ref{alg:enhancedAEM}. When the condition on line \ref{line:nupdate} of Algorithm \ref{alg:enhancedAEM} is met, these quantities will be updated in
sequence to produce $\{(\overline{v}_i, \overline{w}_i)\}_{i=1}^{p}$
using the most recently computed quantities.  In particular, suppose the
updated pairs $\{(\overline{v}_i, \overline{w}_i)\}_{i=1}^{l-1}$ have been
computed. Then the $l$th pair $(v_l, w_l)$ is updated as follows.  First, given
$w_l$, the update  $\overline{v}_l$ is computed by solving 
\begin{align}
\mathcal A_w(w_l, w_l) \overline{v}_l = B w_l - {\sum_{i=1}^{l-1} \mathcal A_w( w_l,  \overline{w}_i) \overline{v}_i} - \sum_{i=l+1}^{p} \mathcal A_w( w_l, w_i) v_i.\label{eq:vp_gsupdate}
\end{align}
Then given $\overline{v}_l$, $\overline{w}_l$ is computed by solving
\begin{align}
\mathcal A_v(\overline{v}_l,\overline{v}_l) \overline{w}_l = B\tran \overline{v}_l 
- \sum_{i=1}^{l-1} \mathcal A_v( \overline{v}_l, \overline{v}_i) \overline{w}_i 
- \sum_{i=l+1}^{p} \mathcal A_v( \overline{v}_l, {v}_i) {w}_i.\label{eq:wp_gsupdate}
\end{align}

\begin{algorithm}[!t]
\caption{\GSnameAbbrv\ enhancement} 
\hspace*{\algorithmicindent} \textbf{Input:} $V_p$ and $W_p$
\begin{algorithmic}[1]
%\Function{GSE}{$V_p,W_p$}
\Function{\GSnameAbbrv}{$V_p,W_p$}
  \For{$l=1,\ldots,p$}
    \State $\overline {v}_l \leftarrow$ solution of equation \eqref{eq:vp_gsupdate} 
    \State $\overline{w}_l\leftarrow$ solution of equation \eqref{eq:wp_gsupdate} 
  \EndFor
  \State $V_p \leftarrow \overline{V}_p$, $W_p \leftarrow \overline{W}_p$
\EndFunction
\end{algorithmic}
\label{alg:GS_AEM}
\end{algorithm}

With $p=2$ as an example, in  vector format, the first step of this %the GS-like
enhancement is to update $(v_1, w_1)$ to $(\overline{v}_1,\overline{w}_1)$ by
solving the following equations: 
\begin{align*}
\left[
\begin{array}{*4{C{2.4cm}}}
\!\!\!A_w({w}_1,{w}_1)\!\!\! & \!\!\!A_w({w}_1,w_2)\!\!\!\\ %& \!\!\!A_w({w}_1,w_3)\\
\phantom{\!\!\!A_w({w}_1,{w}_1)\!\!\!}\\
%\phantom{\!\!\!A_w({w}_1,{w}_1)\!\!\!}
\end{array}
\right]\!\!\!
\left[
\begin{array}{*4{C{.35cm}}}
\!\!\overline{v}_1 \!\!\!\\
\!\!{v}_2 \!\!\!%\\
%\!\!{v}_3 \!\!\!
\end{array}
\right]
&\!\!= \!\!
\left[
\begin{array}{*4{C{1.1cm}}}
\!\!\!B w_1\!\!\!\\
\phantom{\!\!\!B w_1\!\!\!}%\\
%\phantom{\!\!\!B w_1\!\!\!}
\end{array}
\right], \nonumber
\\
\left[
\begin{array}{*4{C{2.4cm}}}
\!\!A_v(\overline{v}_1,\overline{v}_1) \!\!&\!\! A_v(\overline{v}_1,v_2) \!\!\\%&\!\! A_v(\overline{v}_1,v_3)\!\!\\
\phantom{\!\!A_v(\overline{v}_1,\overline{v}_1) \!\!}\\
%\phantom{\!\!A_v(\overline{v}_1,\overline{v}_1) \!\!}
\end{array}
\right]\!\!\!
\left[
\begin{array}{*4{C{.35cm}}}
\!\!\!\overline{w}_1\!\!\! \\
\!\!\!w_2\!\!\! %\\
%\!\!\!w_3\!\!\!
\end{array}
\right]
&\!\!= \!\!
\left[
\begin{array}{*4{C{1.1cm}}}
\!\!\!B\tran \overline{v}_1\!\!\!\\
\phantom{\!\!\!B\tran \overline{v}_1\!\!\!}%\\
%\phantom{\!\!\!B\tran \overline{v}_1\!\!\!}
\end{array}
\right],
\end{align*}
and the second step is to update $(v_2, w_2)$ to
$(\overline{v}_2,\overline{w}_2)$ by solving the second row of  the following equations:  
\begin{align*}
\begin{split}
\left[
\begin{array}{*4{C{2.4cm}}}
\!\!\!A_w(\overline{w}_1,\overline{w}_1)\!\!\! & \!\!\!A_w(\overline{w}_1,{w}_2)\!\!\! \\%& \!\!\!A_w(\overline{w}_1,w_3)\\
\!\!\!A_w(w_2,\overline{w}_1)\!\!\! &\!\!\! A_w(w_2,w_2)\!\!\! \\%&\!\!\!A_w(w_2,w_3)\!\!\!\\
%\phantom{\!\!\!A_w(w_2,\overline{w}_1)\!\!\!}
\end{array}
\right]\!\!\!
\left[
\begin{array}{*4{C{.35cm}}}
\!\!\overline{v}_1 \!\!\!\\
\!\!\overline{v}_2 \!\!\!%\\
%\!\!v_3\!\!\!
\end{array}
\right]
&\!\!= \!\!
\left[
\begin{array}{*4{C{1.1cm}}}
\!\!\!B \overline{w}_1\!\!\!\\
\!\!\!B w_2\!\!\!%\\
%\phantom{\!\!\!B w_2\!\!\!}
\end{array}
\right], 
\\
\left[
\begin{array}{*4{C{2.4cm}}}
\!\!A_v(\overline{v}_1,\overline{v}_1) \!\!&\!\! A_v(\overline{v}_1,\overline{v}_2) \!\!\\%&\!\! A_v(\overline{v}_1,{v}_3)\!\!\\
\!\!A_v(\overline{v}_2,\overline{v}_1) \!\!&\!\! A_v(\overline{v}_2,\overline{v}_2) \!\!%&\!\! A_v(\overline{v}_2,v_3)\!\!\\
%\phantom{\!\!A_v(\overline{v}_2,\overline{v}_1) \!\!}
\end{array}
\right]\!\!\!
\left[
\begin{array}{*4{C{.35cm}}}
\!\!\!\overline{w}_1\!\!\! \\
\!\!\!\overline{w}_2\!\!\! %\\
%\!\!\!w_3\!\!\!
\end{array}
\right]
&\!\!= \!\!
\left[
\begin{array}{*4{C{1.1cm}}}
B\tran \overline{v}_1\\
B\tran \overline{v}_2%\\
%\phantom{B\tran \overline{v}_2}
\end{array}
\right]. 
\end{split}
\end{align*}

This strategy, which we call the \GSnameAbbrv\ enhancement, is summarized in Algorithm \ref{alg:GS_AEM}.
It is an alternative to Algorithm \ref{alg:aem_sr} and is also applied every
 $n_{\text{update}}$ outer iterations.
For a comparison of  Algorithms \ref{alg:aem_sr} and \ref{alg:GS_AEM}, note that 
 Algorithm \ref{alg:aem_sr} (\pgdupdatename) works with a larger system but it can exploit the 
 matricized representation \eqref{eq:wq_Aspd_os}.
Once  the system matrices $\widetilde G_i  = \widetilde W_p\tran G_i \widetilde
W_p$ or $\widetilde K_i= \widetilde V_p\tran K_i \widetilde V_p$ are
formed, if it is not too large, the system in \eqref{eq:wq_Aspd_os} (of order $n_2p$ in this example) 
can be solved using a 
single application of an iterative method such as the preconditioned conjugate gradient (PCG) method.
In contrast, Algorithm \ref{alg:GS_AEM} (\GSnameAbbrv) requires sequential updates of individual components
in equations \eqref{eq:vp_gsupdate}-\eqref{eq:wp_gsupdate}, but with smaller blocks, of order $n_1$ and $n_2$.
As we will show in Section \ref{sec:numexp}, the \GSaemnameAbbrv\ AEM method exhibits  better 
performance in some error measures.

We have found that in practice, the enhancement procedure can be improved by updating 
only a chosen subset of solution pairs rather than all the solution pairs
$\{(v_i,w_i)\}_{i=1}^{p}$.
We discuss a criterion to choose such a subset next.

\subsection{\RstagenameCap\ AEM method} \label{sec:adaptive}
The third enhancement procedure
excerpts and modifies certain computations in  the \stagenameCap\, AEM
method (Lines \ref{ll:Vp} and \ref{ll:Wp}  in Algorithm \ref{alg:stage_aem})
in a computationally efficient way. The procedure
adaptively chooses solution pairs to be updated and solves reduced systems to
update only those pairs.  Let us assume for now that a subset of the solution
pairs to be updated has been chosen. Denote the set of indices of those solution
pairs  by $\ell(p) \subseteq \{1,\ldots,p-1\}$ and the remaining indices by
$\ell^{\text{c}}(p) = \{1,\ldots,p-1\} \setminus \ell(p)$. Then the update is
performed by solving the following equations for
$\overline{V}_{\ell(p)}$ and $\overline{W}_{\ell(p)}$:
\begin{align}\label{eq:vp_as_update}
\sum_{i=0}^{m} (K_i) \overline{V}_{\ell(p)} (\widetilde W_{\ell(p)}\tran G_i  \widetilde W_{\ell(p)})\tran &= B \widetilde W_{\ell(p)} - \sum_{i=0}^m (K_i) V_{\ell^{\text{c}}(p)} (\widetilde W_{\ell(p)}\tran G_i   W_{\ell^{\text{c}}(p)})\tran,
\end{align}
where $\widetilde W_{\ell(p)}$ is obtained by orthonormalizing the columns of $W_{\ell(p)}$, and
\begin{align}\label{eq:wp_as_update}
\sum_{i=0}^{m} (\widetilde V_{\ell(p)}\tran K_i \widetilde V_{\ell(p)}) \overline{W}_{\ell(p)}\tran ( G_i  \tran) &= \widetilde V_{\ell(p)}\tran B  - \sum_{i=0}^m (\widetilde V_{\ell(p)}\tran K_i  V_{\ell^{\text{c}}(p)}) W_{\ell^{\text{c}}(p)}\tran ( G_i \tran),
\end{align}
where $\widetilde V_{\ell(p)}$ is obtained by orthonormalizing the columns of
$\overline{V}_{\ell(p)}$. Then, ${V}_{\ell(p)}$ and ${W}_{\ell(p)}$ are updated
to $\overline{V}_{\ell(p)}$ and $\overline{W}_{\ell(p)}$, while
${V}_{\ell^{\text{c}}(p)}$ and ${W}_{\ell^{\text{c}}(p)}$ remain the same.

We now describe a criterion to choose a subset of the solution pairs to be
updated. Let us assume that $p-1$ iterations of Algorithm \ref{alg:enhancedAEM}
have been performed, and $V_{p-1}$ and $W_{p-1}$ have been computed. The $p$th solution pair ($v_p$, $w_p$) is then computed via
Algorithm \ref{alg:rank1_correction}. If  $p \mod n_{\text{update}}=0$, then a subset of the previous $p-1$
solution pairs is chosen by inspecting the angles between $v_p$ and the columns of $V_{p-1}$ and similarly for $w_{p}$ and $W_{p-1}$. We normalize all vectors $\tilde {v}_i = \frac{v_i}{\|v_i\|_2}$ and compute $\beta_V =  \widetilde V_{p-1}\tran \tilde v_p \in \mathbb{R}^{p-1}$ (the vector of cosines of the angles), and an analogous vector $\beta_W$ using $w_p$ and $W_{p-1}$. The entries of $\beta_V$ and $\beta_W$ indicate how far from orthogonal all previous vectors are to $v_p$ and $w_p$.  Ideally, we want the method to compute $p$ left and right singular vectors of the solution $U$ %in which case $\beta_V = \beta_W=0$.  
(i.e., $\!\beta_V\! =\! \beta_W\!=\!0$).
As the aim is to find good basis vectors for approximating $U$, it is undesirable to keep vectors that are far from being orthogonal to $v_p$ and $w_p$. To resolve this, we choose a subset
of columns of $V_{p-1}$ and $W_{p-1}$ for which the entries of $\beta_V$ and $\beta_W$ are too large; we fix $\tau>0$ and choose
\begin{equation*}%\label{eq:ellp_criteria}
\ell(p) = \{ i \in \{1,\ldots, p-1 \} \mid | [\beta_V]_i | > \tau \text{ or } | [\beta_W]_i | > \tau \}. 
\end{equation*}
Algorithm \ref{alg:astage_AEM} summarizes the resulting
\Rstagename\ (\RstagenameAbbrv) enhancement.
%adaptive-stagewise enhancement procedure.
\begin{algorithm}[!t]
\caption{\RstagenameCap\ enhancement}
\hspace*{\algorithmicindent} \textbf{Input:} $V_p$, $W_p$, and $\tau$
\begin{algorithmic}[1]
\Function{Rstagep}{$V_p,W_p,\tau$}
  \State {Normalize the columns: $\tilde v_{i} = \frac{ v_{i}}{\| v_i\|_2} $, $\tilde w_{i} =  \frac{ w_{i}}{\| w_i\|_2} $ for $i=1,\ldots,p$} \vspace{1mm}
  \State Compute $\beta_V =  \widetilde V_{p-1}\tran \tilde v_p $, $\beta_W =  \widetilde W_{p-1}\tran \tilde w_p$ \vspace{1mm}
  \State Select $\ell(p) = \{ i \in [1,\ldots, p-1] \mid | [\beta_V]_i | > \tau \text{ or } | [\beta_W]_i | > \tau \}$ \vspace{1mm}
  \State $\widetilde W_{\ell(p)} \leftarrow $ orthonormalize $W_{\ell(p)}$ 
  \State $\overline{V}_{\ell(p)} \leftarrow$ solve \eqref{eq:vp_as_update}
  \State $\widetilde V_{\ell(p)} \leftarrow$ orthonormalize $\overline{V}_{\ell(p)}$
  \State $\overline{W}_{\ell(p)} \leftarrow$  solve \eqref{eq:wp_as_update}
  \State $V_{\ell(p)} = \overline{V}_{\ell(p)}$, $W_{\ell(p)} = \overline{W}_{\ell(p)}$
\EndFunction
\end{algorithmic}
\label{alg:astage_AEM}
\end{algorithm}

\section{Numerical experiments}\label{sec:numexp}
In this section, we present the results of numerical experiments with the
algorithms described in Sections \ref{sec:aem} and \ref{sec:update}.  For
benchmark problems, we consider stochastic diffusion problems, where the stochasticity is assumed to be characterized by a prescribed set of real-valued random variables.
We apply suitable stochastic Galerkin finite element discretizations to these problems, which results in linear multi-term matrix equations of the form
\eqref{eq:mat_sys} whose system matrices are symmetric positive-definite.  All
numerical experiments are performed on an \textsc{Intel} 3.1 GHz i7 CPU, with 16 GB
RAM, using \textsc{Matlab} R2019b.

\subsection{Stochastic Diffusion Problems}
Let $(\Omega,\mathcal F, P)$ be a probability space and let $D=[0, 1]\times [0,1]$ be the spatial domain. Next, let $\xi_{i}: \Omega \to \Gamma_{i} \subset \mathbb{R}$, for $i=1,\ldots, m,$ be independent and identically distributed random variables and define $\xi=[\xi_{1}, \ldots, \xi_{m}]$. Then, $\xi: \Omega \to \Gamma$ where $\Gamma \equiv \prod_{i=1}^{m} \Gamma_i$ denotes the image. Given a second-order random field $a: D\times \Gamma \to \mathbb{R}$, we consider the following boundary value problem with constant forcing term $f(x)=1$. Find $u: D \times \Gamma \to \mathbb{R}$ such that
\begin{equation}\label{eq:str_diff}
\left\{
\begin{array}{r l l}
-\nabla \cdot (a(x,\xi) \nabla u(x,\xi) ) &= f(x) &\text{ in } D \times \Gamma,\\
u(x,\xi) &= 0 &\text{ on } \partial D \times \Gamma.
\end{array}
\right.
\end{equation} 
In particular, we will assume that the input random field $a(x,\xi)$ has the affine form
\begin{equation} \label{affine}
a(x,\xi) =  a_{0}(x) +   \sum_{i=1}^{m}  a_i(x) \xi_i,
\end{equation} 
which has the same structure as a truncated Karhunen-Lo\`eve (KL) expansion~\cite{loeve1978probability}, and we will choose the $\xi_i$ to be independent uniform random variables. Recall that if we denote the joint probability density function of $\xi$ by $\rho(\xi)$ then the expected value of a random function $v(\xi)$ on $\Gamma$ is $\langle v \rangle_\rho = \int_\Gamma v(\xi) \rho(\xi) d\xi.$

For the discretization, we consider the stochastic Galerkin method
\cite{babuska2004galerkin,ghanem2003stochastic,lord2014introduction,
xiu2010numerical}, which seeks an approximation to the solution of the following weak
formulation of \eqref{eq:str_diff}: Find  $u(x,\xi)$ in  $V = H_0^1(D) \otimes
L_{\rho}^2(\Gamma)$ such that 
\begin{equation}\label{eq:weak_diff} \left\langle \int_D
a({x},\,\xi) \nabla u({x},\,\xi) \cdot \nabla v({x},\,\xi) d{x}
\right\rangle_\rho = \left\langle\int_D f(x) v({x},\,\xi) dx
\right\rangle_\rho\!\!,  \quad  \forall v \in V.  
\end{equation} 
In
particular, we seek a finite-dimensional approximation of the solution of the
form 
%\begin{equation}\label{eq:sg_sol} 
$\tilde{u}({x},\,{\xi}) = \sum_{s=1}^{n_\xi}
\sum_{r=1}^{n_x} u_{rs} \phi_r({x}) \psi_s({\xi}), $
%\end{equation} 
where
$\{\phi_r\}_{r=1}^{n_x}$ is a set of standard finite element basis functions,
which arises from using continuous piecewise bilinear approximation on a
uniform mesh of square elements (Q1 elements\footnote{Our implementation uses the Incompressible Flow \& Iterative Solver Software (IFISS) \cite{elman2014ifiss,ifiss}.}) and $n_x$ is related to the refinement level of the
spatial mesh. In addition, $\{ \psi_s \}_{s=1}^{n_\xi}$ is chosen to be a finite subset of the set of orthonormal
polynomials that provides a basis for $L_{\rho}^{2}(\Gamma)$ (also known as a generalized
polynomial chaos (gPC), \cite{xiu2002wiener}).  
As the random variables are uniformly distributed, we use
$m$-variate normalized Legendre polynomials $\{\psi_s \}_{s=1}^{n_\xi}$, which are constructed as products of univariate Legendre polynomials, $\psi_s(\xi) =
\prod_{i=1}^{m} \pi_{d_i(s)}(\xi_i)$. Here, $d(s) = (d_1(s),\ldots,d_m(s))$ is
a multi-index and $\pi_{d_i(s)}$ is the $d_i(s)$-order
univariate Legendre polynomial in $\xi_i$. 
A set of 
multi-indices $\{d(s)\}_{s=1}^{n_\xi}$ is specified as a set 
$\Lambda_{m,\, d_\text{tot}} = \{d(s) \in \mathbb{N}^{m}_0: \|{d}(s)\|_1 \leq
d_\text{tot}\}$,  where $\mathbb{N}_0$ is the set of non-negative integers,
$\|{d}(s)\|_1 = \sum_{j=1}^{m} d_j(s)$, and $d_\text{tot}$ defines the maximal
degree of $\{\psi_s(\xi)\}_{s=1}^{n_\xi}$. With this setting, the number of gPC
basis functions is $n_\xi = \text{dim}(\Lambda_{m,\, d_\text{tot}}) =
\frac{(m+d_\text{tot})!}{m!d_\text{tot}!}$.  

Employing a Galerkin projection to \eqref{eq:weak_diff} onto the chosen finite-dimensional space (i.e., using the same
test basis functions as the trial basis functions) and ordering the
coefficients of %\eqref{eq:sg_sol} 
the solution expansion as $u = [u_{11}, \ldots, u_{n_x1},u_{12}, \ldots, u_{n_x n_\xi}]\tran$
results in 
\begin{equation}
\label{eq:linsys_tensor}
\left(\sum_{i=0}^m G_i \otimes K_i \right) u = g_0 \otimes f_0,
\end{equation}
where the system matrices are defined as
\begin{alignat*}{3}
%\label{eq:stiff_matrices}
%\begin{split}
[G_0]_{st} &= \left\langle   \psi_s ({\xi})  \psi_t ({\xi}) \right\rangle_\rho, &&\quad [K_0]_{k\ell} &&= \int_D a_{0}(x) \nabla \phi_k ({x}) \cdot \nabla \phi_\ell ({x}) d{x}, \\
 [G_i]_{st} &=  \left\langle  \xi_i\, \psi_s ({\xi})  \psi_t ({\xi}) \right\rangle_\rho, &&\quad [K_i]_{k\ell} &&= \int_D  a_i({x}) \nabla \phi_k ({x}) \cdot  \nabla \phi_\ell ({x}) d{x}, 
%\end{split}
\end{alignat*} 
for $i = 1,\,\ldots,\,m$, $s,t=1,\ldots, n_{\xi}$ and $k,\ell=1,\ldots n_{x}$. Due to the deterministic forcing term $f(x)=1$,
the right-hand side has a rank-one structure (i.e., $r=0$ in 
\eqref{eq:kron_sys}), with $[f_0]_k = \int_D f(x) \phi_k({x}) d{x},$ and  $[g_0]_s = \left\langle  \psi_s ({\xi}) \right\rangle_\rho$.
Matricizing
\eqref{eq:linsys_tensor} gives the multi-term matrix equation as shown in
\eqref{eq:mat_sys} with $n_1 = n_x$ and $n_2 = n_\xi$, and now we can apply the
AEM methods to compute an approximate solution of the equation.

\subsection{Benchmark problem 1: separable exponential covariance}  \label{sec:ex_spd}

In this problem, we assume that the random field $a(x,\xi)$ is a truncated KL expansion
\begin{equation} \label{eq:def_rf}
a(x,\xi) =  \mu + \sigma \sum_{i=1}^{m}  \sqrt{\lambda_i} \varphi_i(x) \xi_i,
\end{equation} 
where $\mu$ is the mean of $a(x,\xi)$, $\{(\varphi_i(x),\lambda_i)\}_{i=1}^{m}$ are eigenpairs of the integral operator
associated with the separable covariance kernel $C({x},\, {y}) \equiv  \exp \left
( - \frac{|x_1 - y_1|}{c} - \frac{|x_2 - y_2|}{c} \right)$, $c$ is the associated correlation
length, and $\sigma^2$ is the variance of the untruncated random field. In addition, each $\xi_i \sim U(-\sqrt{3},\sqrt{3})$ and so has mean zero and variance one.

In the following sections, we compare the five AEM variants, 
\stagenameCap\ (Algorithm \ref{alg:stage_aem}), 
\baselinename\ (Algorithm \ref{alg:base_aem}), 
\pgdaemname\  (Algorithm \ref{alg:aem_sr}),  
\GSaemnameAbbrv\  (Algorithm \ref{alg:GS_AEM}),  
and
\Rstagename\  (Algorithm \ref{alg:astage_AEM}).  For orthonormalization in \pgdaemname\  (Algorithm \ref{alg:aem_sr}) and \Rstagename\  (Algorithm \ref{alg:astage_AEM}), we use \textsc{Matlab}'s \texttt{qr} function. For assessing performances, we explore two key aspects. The first is the accuracy of the computed solutions,
which we assess by computing two error metrics: cosines of angles
between the truth singular vectors and the columns of the computed
factors (Section \ref{sec:svs}), and errors between the
 truth solution and the computed solution measured in three different norms (Section
\ref{sec:acc}). The second aspect is timings and scalability (Section
\ref{sec:time}). As the assessment of the first aspect requires the ground
truth solution of \eqref{eq:linsys_tensor}, which is computed using
\textsc{Matlab}'s backslash operator, and its singular vectors, we choose
small-sized problems in Sections \ref{sec:svs}--\ref{sec:acc}. When making comparisons with the truth solution, we set the maximum number of outer iterations for all the AEM methods to be $p_{\max}=\min(n_x,n_\xi)=56$. Larger problems are considered in Section \ref{sec:time}, where 
scalability matters and finding the truth solution is impossible with the available resources.

\begin{figure}[!h]
%\centering
\begin{minipage}{.325\linewidth}
\centering
\includegraphics[scale=.435]{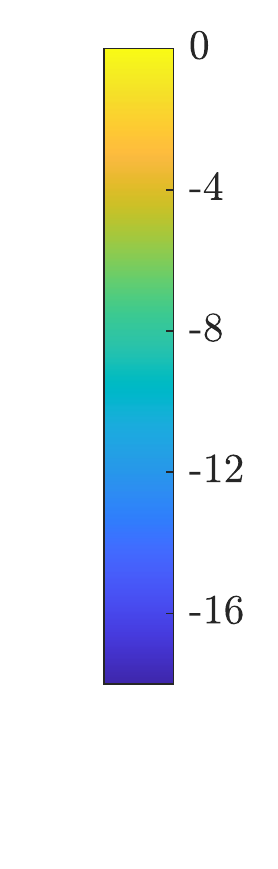}
\end{minipage}
\begin{minipage}{.325\linewidth}
\centering
\subfloat[\stagenameCap, $ V^\ast{}\tran \widetilde V_p$]  {\label{fig:angle_test1a}\includegraphics[scale=.425]{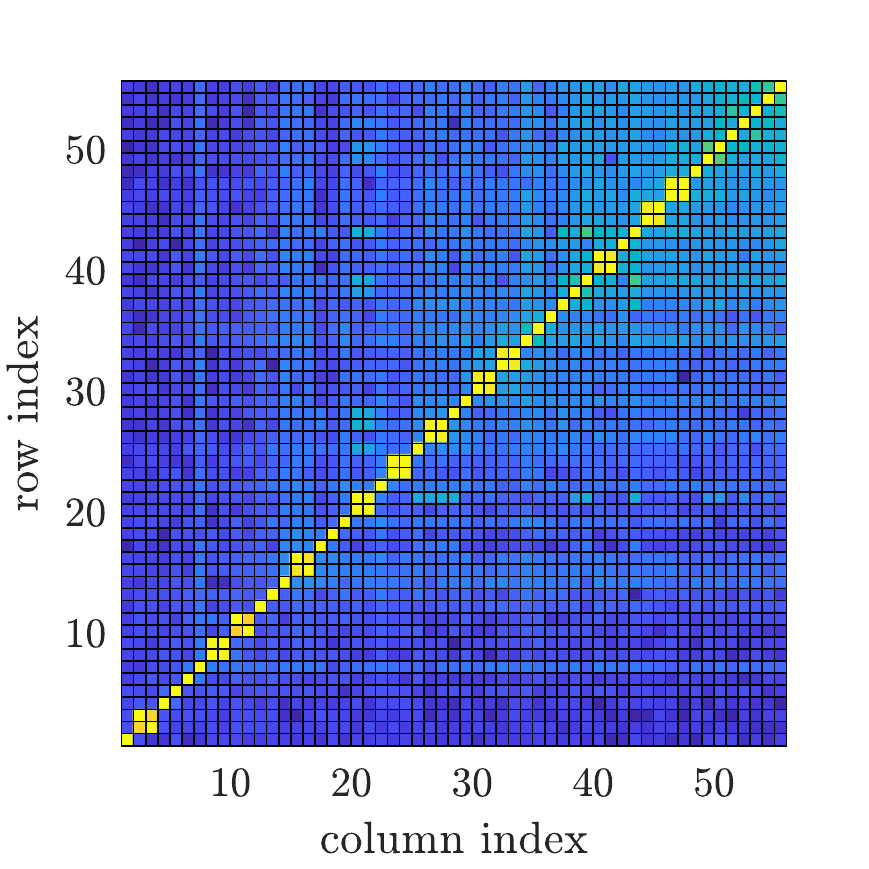}} 
\end{minipage}
\begin{minipage}{.325\linewidth}
\centering
\subfloat[\baselinename, $ V^\ast{}\tran \widetilde V_p$]  {\label{fig:angle_test2a}\includegraphics[scale=.425]{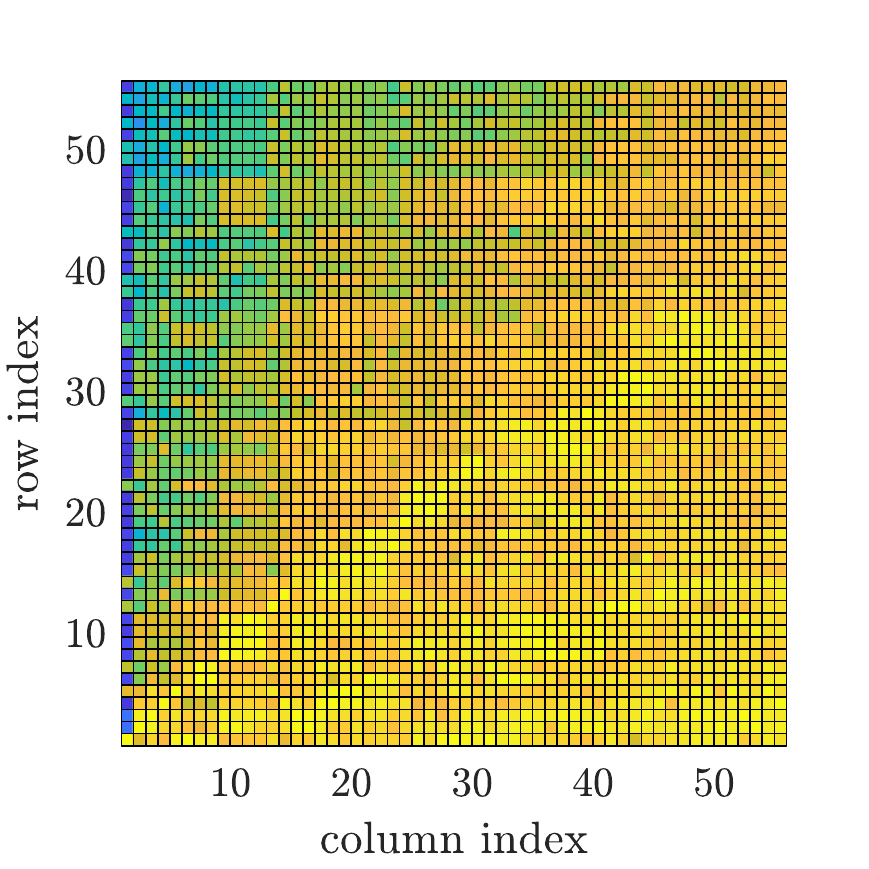}}
\end{minipage}\\
\begin{minipage}{.325\linewidth}
\centering
\subfloat[\pgdupdatename, $V^\ast{}\tran \widetilde V_p$]  {\label{fig:angle_test3a}\includegraphics[scale=.425]{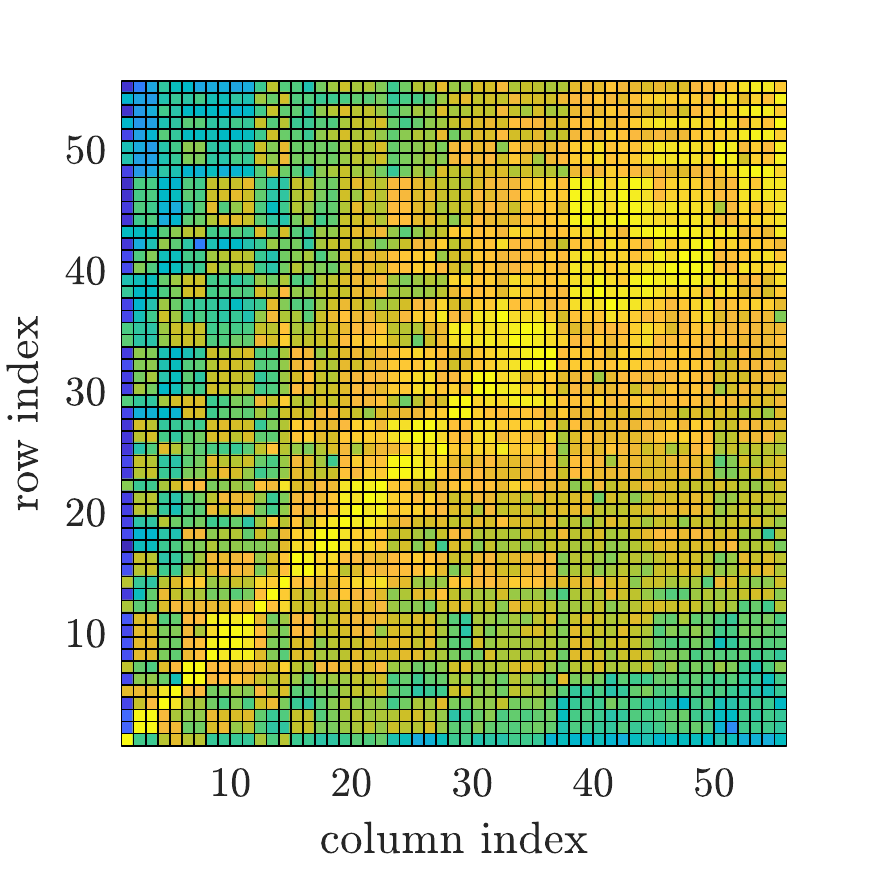}} 
\end{minipage}
\begin{minipage}{.325\linewidth}
\centering
\subfloat[\GSnameAbbrv, $V^\ast{}\tran \widetilde V_p$]  {\label{fig:angle_test4a}\includegraphics[scale=.425]{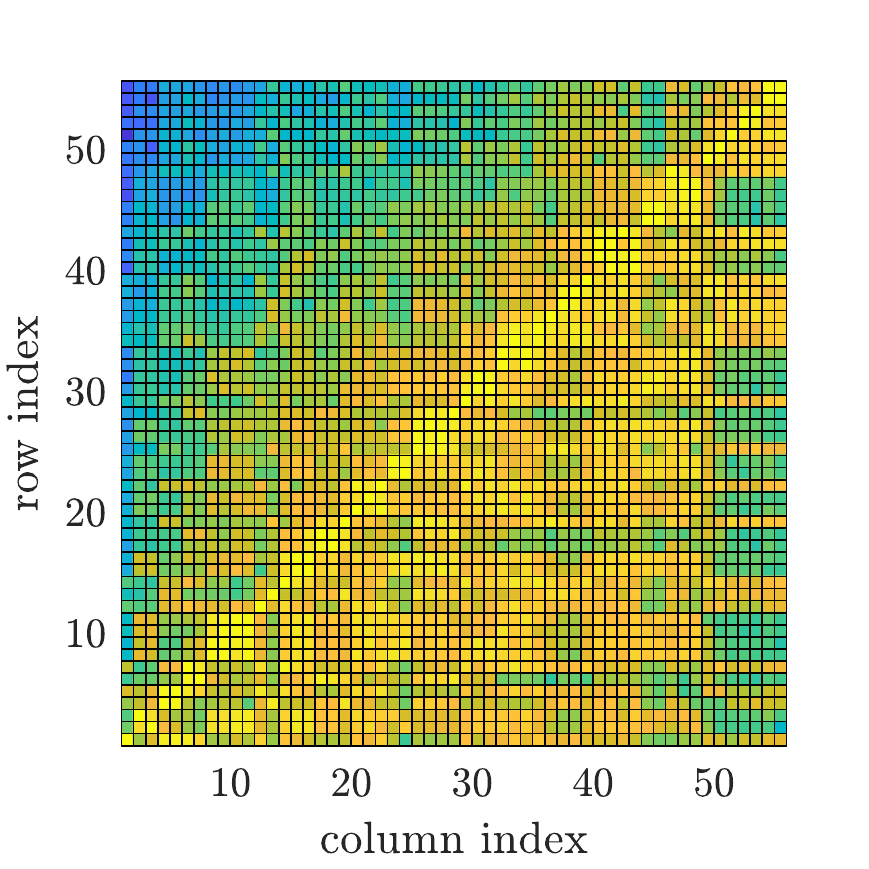}} 
\end{minipage}
\begin{minipage}{.325\linewidth}
\centering
\subfloat[\RstagenameAbbrv, $V^\ast{}\tran \widetilde V_p$]  {\label{fig:angle_test5a}\includegraphics[scale=.425]{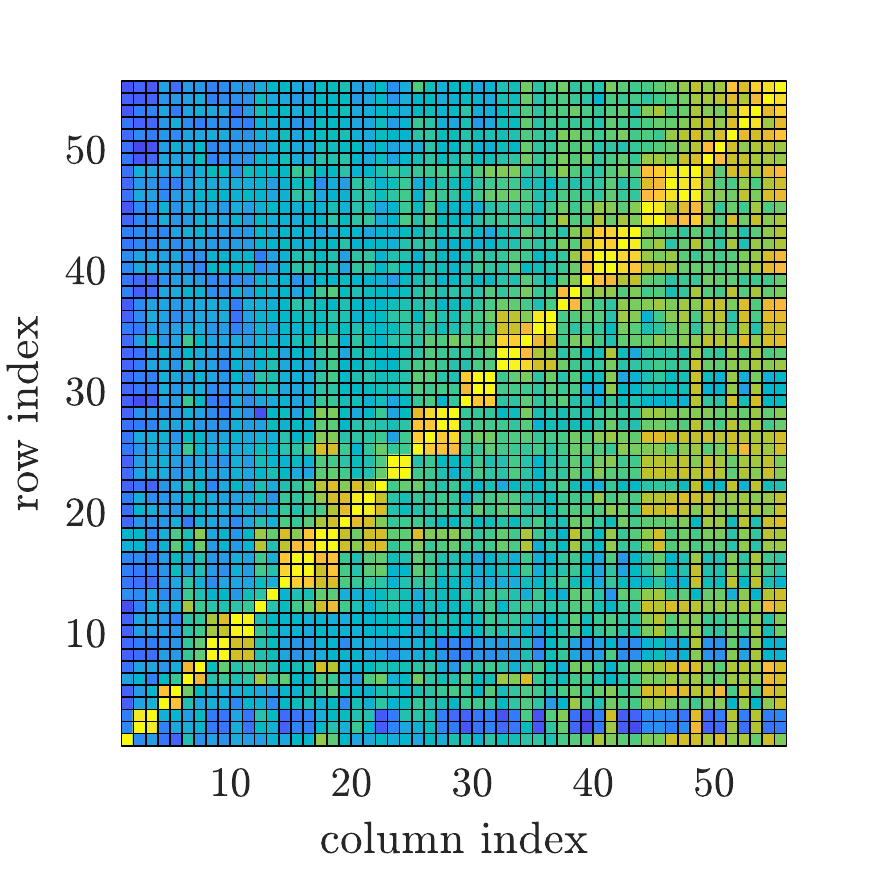}} 
\end{minipage}
\vspace{1mm}
\hrule
\begin{minipage}{.325\linewidth}
\centering
\includegraphics[scale=.435]{figs/colorbar}
\end{minipage}
\begin{minipage}{.325\linewidth}
\centering
\subfloat[\stagenameCap, $ W^\ast{}\tran \widetilde W_p$]  {\label{fig:angle_test1b}\includegraphics[scale=.425]{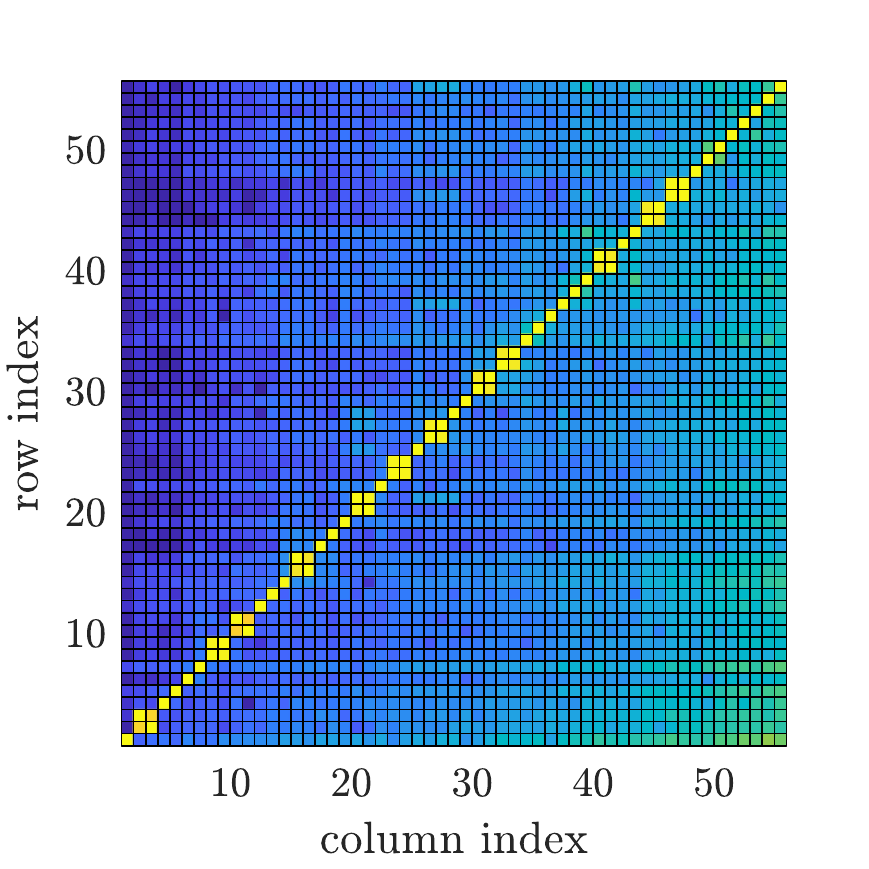}} 
\end{minipage}
\begin{minipage}{.325\linewidth}
\centering
\subfloat[\baselinename, $ W^\ast{}\tran \widetilde W_p$]  {\label{fig:angle_test2b}\includegraphics[scale=.425]{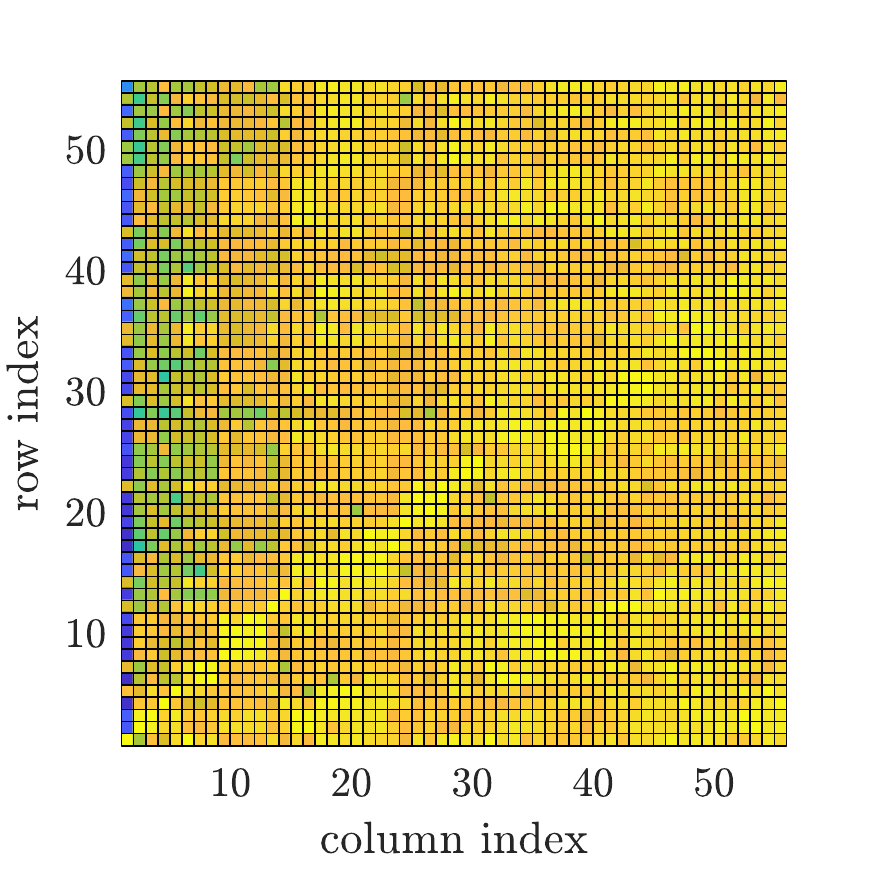}} 
\end{minipage}\\
\begin{minipage}{.325\linewidth}
\centering
\subfloat[\pgdupdatename, $W^\ast{}\tran \widetilde W_p$]  {\label{fig:angle_test3b}\includegraphics[scale=.425]{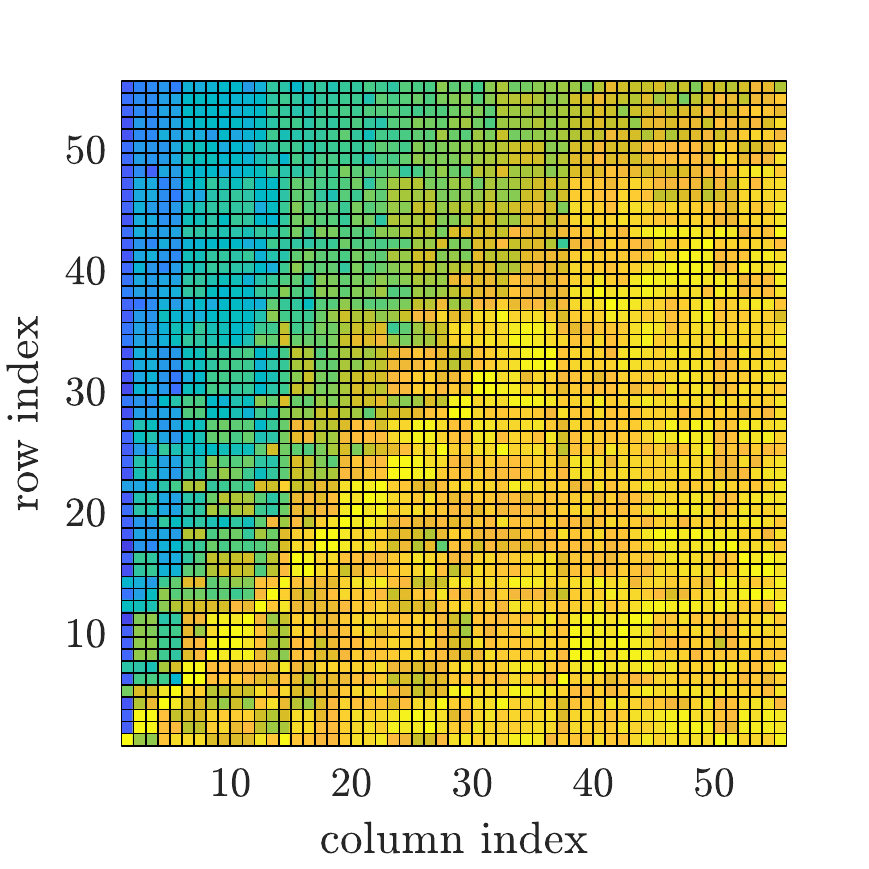}} 
\end{minipage}
\begin{minipage}{.325\linewidth}
\centering
\subfloat[\GSnameAbbrv, $W^\ast{}\tran \widetilde W_p$]  {\label{fig:angle_test4b}\includegraphics[scale=.425]{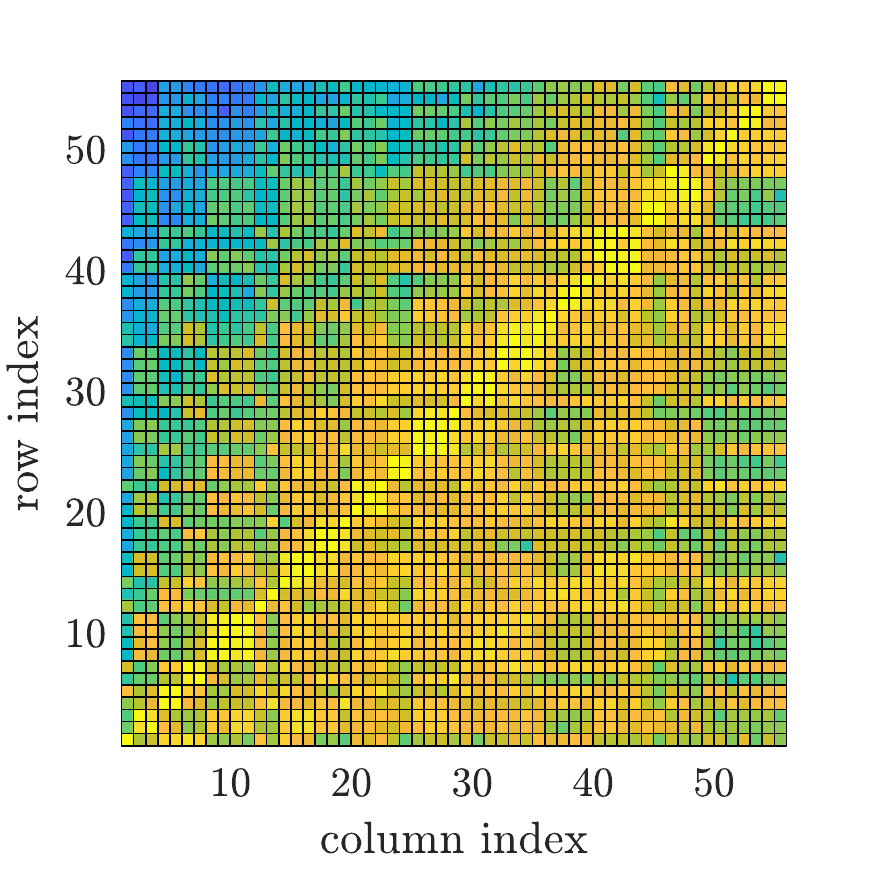}}
\end{minipage}
\begin{minipage}{.325\linewidth}
\centering
\subfloat[\RstagenameAbbrv, $W^\ast{}\tran \widetilde W_p$]  {\label{fig:angle_test5b}\includegraphics[scale=.425]{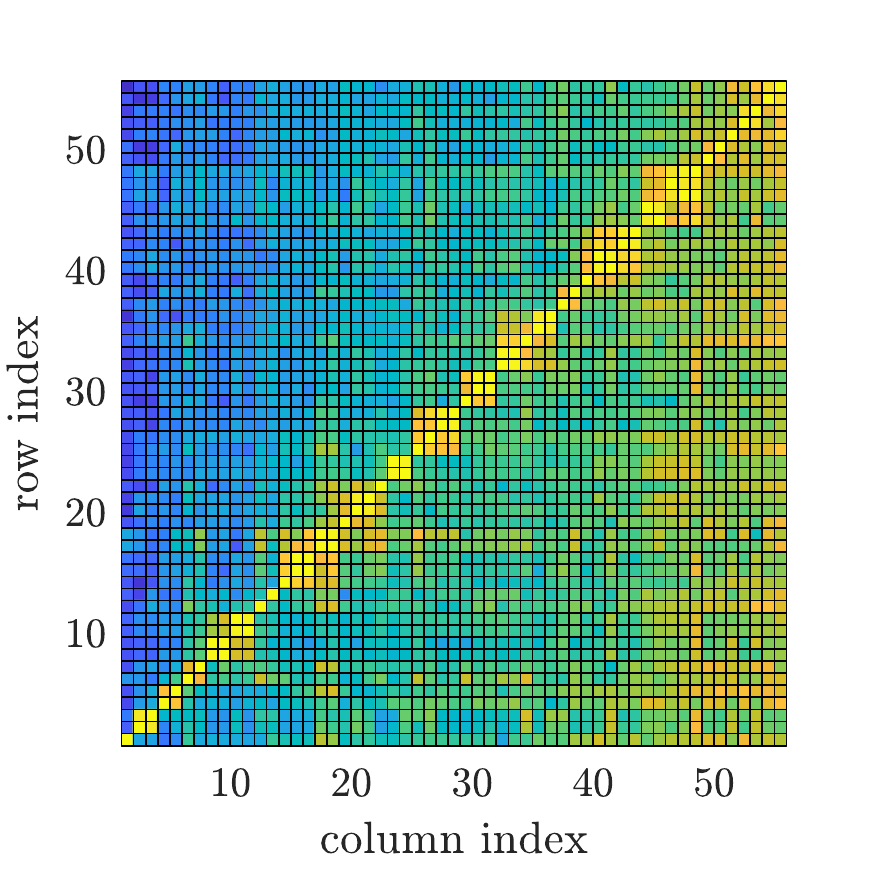}} 
\end{minipage}
\caption{Cosines of angles  (plotted in $\log$ scale) between the left singular vectors $V^\ast$ and $\widetilde V_p$, and the right singular vectors $W^\ast$ and $\widetilde W_p$, where $\widetilde V_p$ and $\widetilde W_p$ are computed using the \stagenameCap\ and \baselinename\ AEM methods, and the \text{EnhancedAEM} methods with  \pgdupdatename, \GSnameAbbrv, and \RstagenameAbbrv\ enhancements.}
\label{fig:angle_test2}
\end{figure}

\subsubsection{Relation to singular vectors} \label{sec:svs}
We begin by exploring how the factors in the approximate solutions
constructed by each of the methods  %under consideration 
compare with the left
and right singular vectors of the true solution matrix $U$. This is
important because (i) singular vectors
represent the most effective choice with respect
to the Frobenius norm for approximating a matrix $U$.
That is, the minimum error over all rank-$p$ approximations is 
$\|U - \widetilde{V}_p\Sigma_p \widetilde{W}_p\tran\|\fro$, where
$U=\widetilde{V}\Sigma \widetilde{W}\tran$ is the singular value
decomposition~\cite{eckart1936approximation}, and (ii) in some applications such as collaborative
filtering for recommendation systems, computing singular vectors accurately is
very important for precise predictions \cite{jain2013low,koren2009matrix,zhou2008large}. 
For these tests,  the diffusion coefficient is given by \eqref{eq:def_rf} with $(\mu,\sigma) = (1,.1)$ and $c=2$.  We use
a spatial discretization with grid level 4 (i.e.,
grid spacing $\frac{1}{2^4}$, and $n_x=225$) and
we truncate the expansion
\eqref{eq:def_rf} at $m=5$. % (i.e., five random variables), which captures
%95.72$\%$ of the total variance of the random field. 
For the stochastic Galerkin
approximation, we choose 
$d_{\text{tot}}=3$ which gives $n_\xi = 56$.

For any approximation of the form \eqref{eq:apprx_rep}, let $\widetilde V_p$ and
$\widetilde W_p$ be normalized versions of the factors, i.e., each column of
$\widetilde V_p$ and $\widetilde W_p$ is scaled to have unit norm. 
From the ground truth solution $U$, the matrices $V^\ast$ and $W^\ast$ of left and right singular vectors are computed.
 The entries of $V^\ast{}\tran \widetilde V_p$, the cosines of the angles between the left singular vectors of the true
solution and the left vectors defining the approximate solution, together with
the analogous angles for the right vectors, $W^{\ast}{}\tran \widetilde W_p$,
give insight into the quality of the approximate solution. Figures
\ref{fig:angle_test1a} and \ref{fig:angle_test1b} and Figures \ref{fig:angle_test2a} and \ref{fig:angle_test2b} depict the cosines of the angles between the
singular vectors and the columns of $\widetilde V_p$ and
$\widetilde W_p$ computed using the \stagenameCap\ AEM and \baselinename\
AEM methods discussed in Section \ref{sec:aem}. 
It can be seen from these results  (in Figures \ref{fig:angle_test1a} and \ref{fig:angle_test1b})  that the \stagenameCap\ AEM 
method does a good job of approximating the singular vectors of the solution. 
That is, the values
of the diagonal entries are close to one and the values of the off-diagonal entries
are close to zero. On the other hand, the \baselinename\ AEM method (see Figures
\ref{fig:angle_test2a} and \ref{fig:angle_test2b}) is far less effective.  
The $2\times 2$ blocks on the diagonals
in Figures \ref{fig:angle_test1a} and \ref{fig:angle_test1b} reflect the
presence of equal singular values.

Figures \ref{fig:angle_test3a}--\ref{fig:angle_test5a} and  \ref{fig:angle_test3b}--\ref{fig:angle_test5b} show analogous results for EnhancedAEM with 
\pgdupdatename\ 
(Algorithm \ref{alg:aem_sr}), 
\GSnameAbbrv\ 
(Algorithm \ref{alg:GS_AEM}), and 
\RstagenameAbbrv\
(Algorithm \ref{alg:astage_AEM}). 
Since we attempt to see each method's best possible results
without considering the computational costs, we set $k_{\max}=5$ and
$n_\text{update}=1$ (i.e., enhancements are performed at every outer iteration) in
Algorithm \ref{alg:enhancedAEM}. For the same reason, we set  \GSnameAbbrv\ 
to update all the
 solution pairs and, for  \RstagenameAbbrv, 
 we set $\tau=.001$.  
 With \pgdupdatename, the spatial %SE, the spatial
component gets reduced (i.e., we form $\widetilde K_i= \widetilde V_p\tran
K_i\widetilde V_p$) and $W_p$ is updated. %The results in 
Figures
\ref{fig:angle_test3a} and \ref{fig:angle_test3b} show that this computation
improves the quality of the resulting factor $W_p$ (and $V_p$ as well) as approximate singular vectors,  compared to those obtained
with the \baselinename\ method. It is
evident that \GSnameAbbrv\ further
improves the quality of $\widetilde V_p$ and $\widetilde W_p$ 
(Figures \ref{fig:angle_test4a} and \ref{fig:angle_test4b}) as
approximate singular vectors, and \RstagenameAbbrv\ 
is nearly as effective as the 
\stagenameCap\ AEM approach (Figures \ref{fig:angle_test5a} and \ref{fig:angle_test5b}).

\subsubsection{Assessment of solution accuracy} \label{sec:acc}
We now compare the convergence behavior of the variants of the AEM methods
introduced in Sections \ref{sec:aem} and \ref{sec:update}. We use two different
settings for the stochastic diffusion coefficient: [exp1] $(\mu,\sigma) = (1,
.1)$, $c=2$ and [exp2] $(\mu,\sigma) = (1, .2)$, $c=.5$. We again truncate
the series \eqref{eq:def_rf} at $m=5$ and, for the Legendre basis polynomials, we consider
$d_{\text{tot}}=3$ which gives $n_\xi = 56$. We deliberately keep
the same value for $m$ and $d_{\text{tot}}$ for both settings so that we can
keep the dimensions of the problem the same and, thus, directly compare the
behavior of each method in different problem settings. We also use the same
parameters for the EnhancedAEM methods as before (i.e., $k_{\max}=5$, $n_\text{update}=1$, and
$\tau=.001$).

\begin{figure}[!t]
%\centering
\begin{minipage}{.475\linewidth}
\centering
\subfloat[Energy norm - Exp1]  {\label{fig:conv_test1a}\includegraphics[scale=.61]{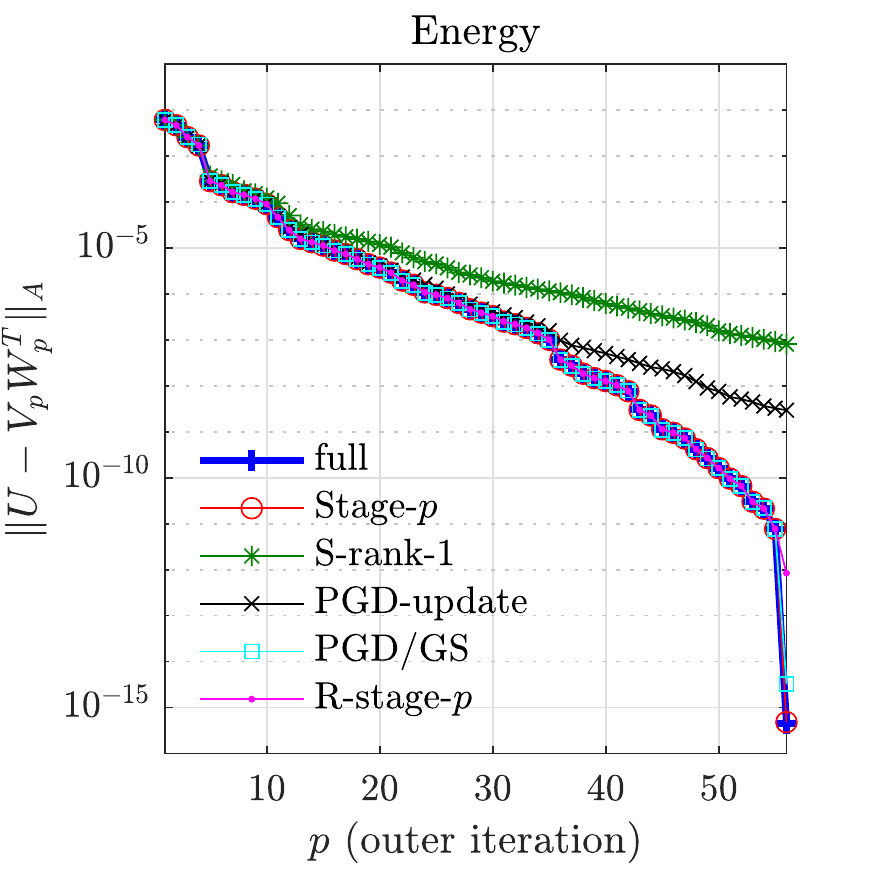}}
%\subfloat[Error norm]  {\label{fig:conv_test1b}\includegraphics[scale=.61]{figs/[EASY]M5_error_vert}}
%\subfloat[Residual norm]  {\label{fig:conv_test1c}\includegraphics[scale=.61]{figs/[EASY]M5_res_vert}} \\
\end{minipage}
%\hspace{5mm}
\begin{minipage}{.475\linewidth}
\centering
\subfloat[Energy norm - Exp2]  {\label{fig:conv_test2a}\includegraphics[scale=.61]{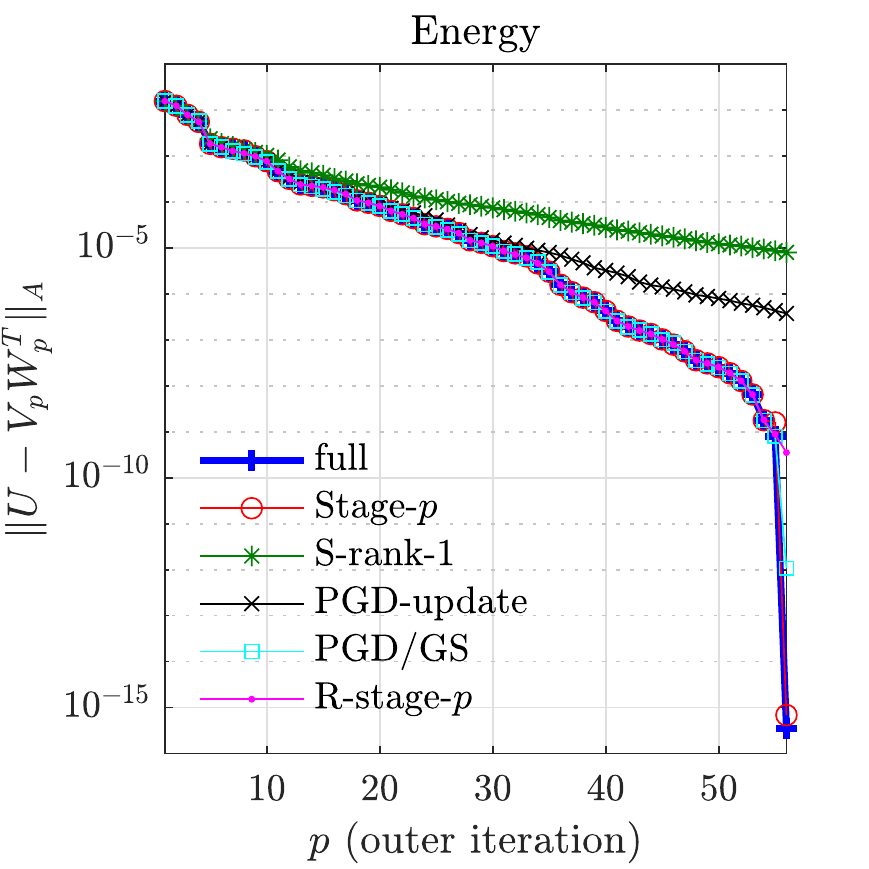}}
\end{minipage}
%\subfloat[Error norm]  {\label{fig:conv_test2b}\includegraphics[scale=.61]{figs/[HARD]M5_error_vert}}
%\subfloat[Residual norm]  {\label{fig:conv_test2c}\includegraphics[scale=.61]{figs/[HARD]M5_res_vert}} 
\caption{Solution errors measured in the energy norm}%different norms: the energy norm%, the error norm and the residual norm.}
\label{fig:conv_test}
\end{figure}

For each method, the approximate solution $\ApproxSol$ is computed and we measure the accuracy compared to the reference solution $U$. We did this using three different
metrics: the energy norm error $\| U - \ApproxSol\|_A$,  the error in the Frobenius norm $\| U - \ApproxSol \|\fro$, and the
residual in the Frobenius norm $\| B - \mathcal A(\ApproxSol)\|\fro$.  
Here, we only report the energy norm errors (in Figure \ref{fig:conv_test}), as behavior for the other two metrics is virtually identical.
For comparison, a rank-$p$ reference solution (referred to as ``full" in Figure \ref{fig:conv_test}) 
is also obtained directly from the first $p$ singular values and singular
vectors of $U$.  

For both settings, as expected, the convergence behavior 
of the \baselinename\ AEM method is significantly worse than that of the rank-$p$ reference solution,
whereas that of the \stagenameCap\ AEM method is virtually the same as for the full direct
solver. The EnhancedAEM method with \pgdupdatename\ 
converges well until a certain level of
accuracy is achieved, but it fails to achieve a high level of accuracy. 
In both experiments, the EnhancedAEM methods with \GSnameAbbrv\ and \RstagenameAbbrv\ 
are more effective than with the \pgdupdatename. The accuracy that those two methods achieve is virtually the same as that of the \stagenameCap\ AEM method and the full direct solver. % in both experimental settings.

\subsubsection{Computational timings} \label{sec:time}
The above results do not account for computational costs; we now investigate timings under various experimental settings.
This is important for large-scale applications, and so we now consider a finer spatial grid,
with grid level 6 (i.e., grid spacing $\frac{1}{2^6}$, and
$n_x=3969$), as well as larger parameter spaces, with %and consider 
$m=\{20,24\}$ (the number of random variables in
\eqref{eq:def_rf}) and $d_{\text{tot}}=4$, which results in $n_\xi =
\{10626,20475\}$. We use the same settings for the stochastic diffusion
coefficient [exp1] $(\mu,\sigma) = (1, .1)$, $c=2$ and [exp2] $(\mu,\sigma) =
(1, .2)$, $c=.5$. Again, we set $m$ and $d_{\text{tot}}$ to be the same for both 
problems, as we want to keep the dimensions fixed  so that we can make direct and fair comparisons. % on the problems with different settings.

Before we present these results, we summarize the
systems of equations to be solved for each of the EnhancedAEM methods
and the adjustable parameters that affect the performances of the
methods.\footnote{The results of using the \stagenameCap\ and  \baselinename\ AEM
methods are not reported because the \stagenameCap\ AEM method is computationally
too expensive and the \baselinename\ AEM method exhibits poor convergence behavior
and, indeed, fails to satisfy the given convergence criterion.} We first describe how
we solve the systems arising at the $p$th outer iteration when
the condition for applying the enhancement %\textsc{Enhancement procedure} 
is met, as well as the systems arising in 
\textsc{RankOneCorrection} (Algorithm \ref{alg:rank1_correction}). We
use PCG to solve each system of equations using mean-based preconditioners
\cite{powell2009block}, which are constructed using reduced versions of the matrices  %only consider the mean filed of the random
$K_0$ and $G_0$, that are adapted to each method. For all
systems, each PCG iteration requires matrix-vector
products in the matricized form (see
\cite{ballani2013projection,kressner2011low,lee2017preconditioned} for detailed
matrix operations)
\begin{equation*}
 \sum_{i=0}^m (M_x^{-1} \widetilde{K}_i) X (M_{\xi}^{-1} \widetilde{G}_i)\tran,
\end{equation*}
where $X$ is a quantity to be updated, $\widetilde{K}_i$ and $\widetilde{G}_i$
are reduced matrices, and $M_x$ and $M_\xi$ are the preconditioner factors. Table
\ref{tab:systems_precond} summarizes each system matrix and
preconditioner.\footnote{Note that, for \pgdupdatename, one can always choose the smallest solution component to update. In practice, however, updating the $W_p$ component (i.e., reduction in $\{K_i\}_{i=0}^{m}$) always requires the smallest computational costs and, thus, we only report the result of updating $W_p$.} 

\begin{table}[tbhp]
{\footnotesize
  \caption{System matrices and preconditioners for each \textsc{Enhancement} procedure}\label{tab:systems_precond}
\begin{center}
\renewcommand{\arraystretch}{1.25}
  \begin{tabular}{|c|c|c|c|c|c|c|} 
\hline
Name & X & $\widetilde K_i$ & $\widetilde G_i$ & $M_x$ & $M_\xi$ & Eqs\\
\hline
\baselinename & $v_p$ & $K_i$ & $ w_p\tran G_i  w_p$ & $K_0$ & 1 & \eqref{eq:vp_rk1}\\
\cline{2-6}
(Alg. \ref{alg:rank1_correction}) &$w_p\tran$ & $ v_p\tran  K_i  v_p$ & $G_i$ &$1$ & $G_0$&  \eqref{eq:wp_rk1} \\
\hline 
\pgdupdatename\ & $V_p$ & $K_i$ & $\widetilde W_p\tran G_i \widetilde W_p$ & $K_0$& $\widetilde W_p\tran G_0 \widetilde W_p$& \\
\cline{2-6}
(Alg. \ref{alg:aem_sr}) &$W_p\tran$ & $\widetilde V_p\tran  K_i \widetilde V_p$ & $G_i$ &$\widetilde V_p\tran  K_0 \widetilde V_p$& $G_0$ &  \eqref{eq:wq_Aspd_os} \\
\hline
\GSnameAbbrv\ & $v_l$ & $K_i$ & $ w_l\tran G_i  w_l$ & $ K_0$& 1 &  \eqref{eq:vp_gsupdate}\\
\cline{2-6}
(Alg. \ref{alg:GS_AEM}) &$w_l\tran$ & $\bar v_l\tran  K_i \bar v_l$ & $G_i$ &1 &$G_0 $&  \eqref{eq:wp_gsupdate} \\
\hline
\RstagenameAbbrv\ & $V_{\ell(p)}$ & $K_i$ & $\!\!\widetilde W_{\ell(p)}\tran G_i \widetilde W_{\ell(p)}\!\!$ & $K_0$& $\!\!\widetilde W_{\ell(p)}\tran G_0 \widetilde W_{\ell(p)}\!\!$&  \eqref{eq:vp_as_update}\\
\cline{2-6}
(Alg. \ref{alg:astage_AEM}) &$W_{\ell(p)}\tran$ & $\widetilde V_{\ell(p)}\tran  K_i \widetilde V_{\ell(p)}$ & $G_i$ &$ \widetilde V_{\ell(p)}\tran  K_0 \widetilde V_{\ell(p)}$& $G_0$ &  \eqref{eq:wp_as_update} \\
\hline
\end{tabular}
\end{center}
}
\end{table}

Now, we discuss adjustable parameters. The EnhancedAEM methods (Algorithms
\ref{alg:enhancedAEM}--\ref{alg:check_conv}) require parameters
$p_{\max}$, $k_{\max}$, $n_{\text{update}}$, and $\epsilon$. We set
$p_{\max}=1000$ to prevent excessive computations. 
We found that
 choosing $k_{\max} > 2$ results in negligible difference in accuracy, but
requires extra computations and, thus, 
we use 
$k_{\max}=\{1,2\}$. %The third parameter 
For $n_{\text{update}}$, which
determines how often the enhancement procedure is called, we vary $n_{\text{update}}$ as $\{5,10,20,30\}$. Next, 
we use $\epsilon$ to check the convergence (as in Algorithm \ref{alg:check_conv}), and
we vary $\epsilon$ as $\{10^{-10},10^{-9},10^{-8},10^{-7}\}$.
Finally, for \GSnameAbbrv\ and \RstagenameAbbrv, we empirically found that choosing $\tau > 0.05$ results
in decreased accuracy in the approximate solution and, thus, we set $\tau =
0.05$.

Next, we set parameters for the PCG method. 
For all systems, the stopping
criterion uses the relative residual in the Frobenius norm. 
We use two different tolerances:
$\tolbasis{}$ for solving systems that arise in \textsc{RankOneCorrection} and \GSnameAbbrv, 
and $\tolblock{}$ for solving systems that arise in \pgdupdatename\ and \RstagenameAbbrv.  
%For
%choosing the values of $\tolbasis{}$ and $\tolblock{}$,  we have tested the
%EnhancedAEM methods for various parameter settings
%$(\tolbasis{},\tolblock{})=\{10^{-8},10^{-7},10^{-6},10^{-5}\}^2$. We have
%observed the following from preliminary numerical experiments: 
We choose the values of $\tolbasis{}$ and $\tolblock{}$ based on results of preliminary numerical experiments with the EnhancedAEM methods for $(\tolbasis{},\tolblock{})=\{10^{-8},10^{-7},10^{-6},10^{-5}\}^2$:
(i)
setting $\tolbasis{} < 10^{-5}$ does not result in improved accuracy of
approximate solutions and, thus we set $\tolbasis{} = 10^{-5}$, %for the
%following experiments, 
and (ii) for a given outer
iteration tolerance $\epsilon$, having too mild PCG tolerance $\tolblock{} >
10^2 \epsilon$ results in poor performance and having stringent tolerance
$\tolblock{} < 10^2 \epsilon$ results in negligible difference in %\marginpar{\hspace{.5in}\HE{OK?}}
accuracy; thus, we use $\tolblock{}=10^2 \epsilon$. % for the following experiments. 
Table \ref{tab:params} summarizes the parameters used for the
experiments.

\begin{table}[tbhp]
{\footnotesize
  \caption{Parameters used in the experiments for measuring timings}\label{tab:params}
\begin{center}
\renewcommand{\arraystretch}{1.25}
  \begin{tabular}{l|l} 
\hline
the maximum number of outer iterations & $p_{\max}=1000$\\
the maximum number of inner iterations & $k_{\max}=\{1,2\}$\\
the frequency of the enhancement procedure & $n_{\text{update}}=\{5,10,20,30\}$\\
the stopping tolerance for outer iterations & $\epsilon \!= \!\{10^{-10},10^{-9},10^{-8},10^{-7}\}\!\!\!\!\!$\\
PCG stopping tolerance for \textsc{RankOneCorrection} and \GSnameAbbrv & $\tolbasis{}=10^{-5}$\\
PCG stopping tolerance for \pgdupdatename\ and \RstagenameAbbrv & $\tolblock{}=10^{2} \epsilon$\\
\hline
\end{tabular}
\end{center}
}
\end{table}

\begin{figure}[!h]
\centering
\includegraphics[scale=.5]{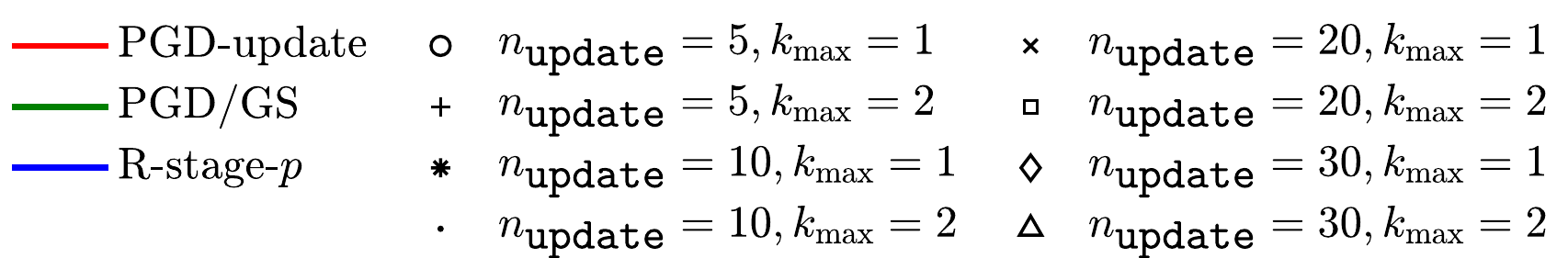} \vspace{-4mm}\\
\subfloat[$m=20$, {[exp1]}]  {\label{fig:timing_easy1b}\includegraphics[scale=.45]{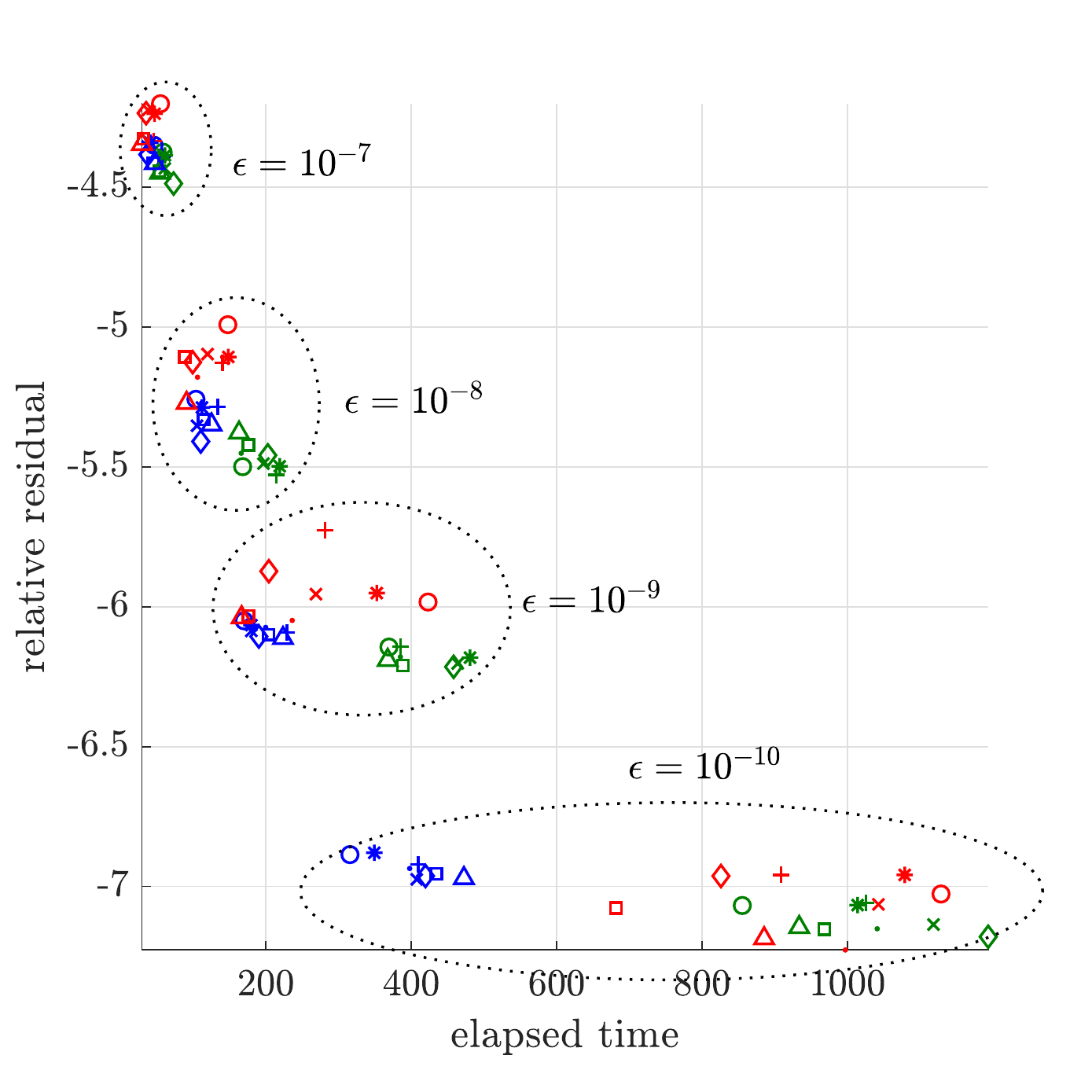}} 
\subfloat[$m=24$, {[exp1]}]  {\label{fig:timing_easy2b}\includegraphics[scale=.45]{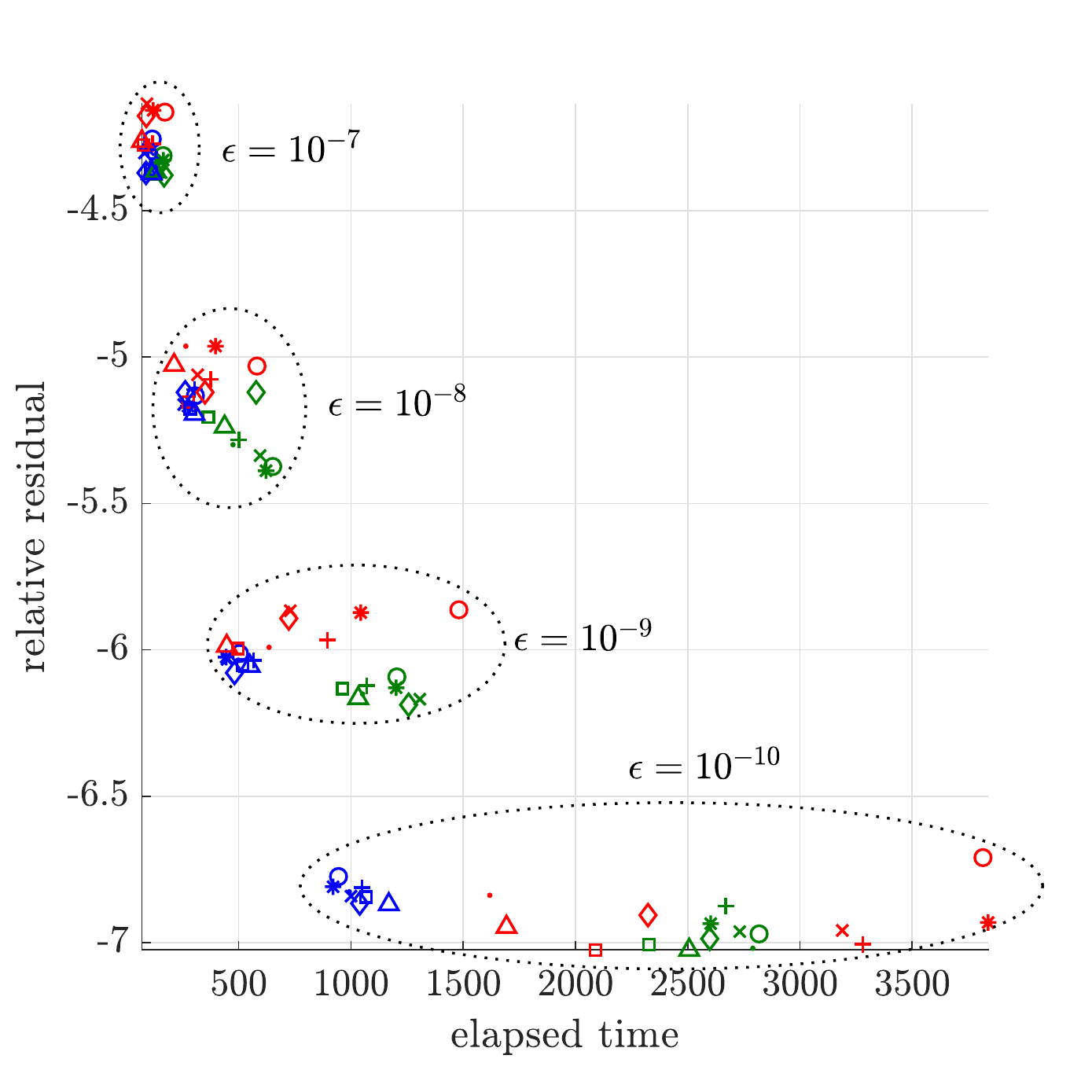}} \\
\subfloat[$m=20$, {[exp2]}]  {\label{fig:timing_hard1b}\includegraphics[scale=.45]{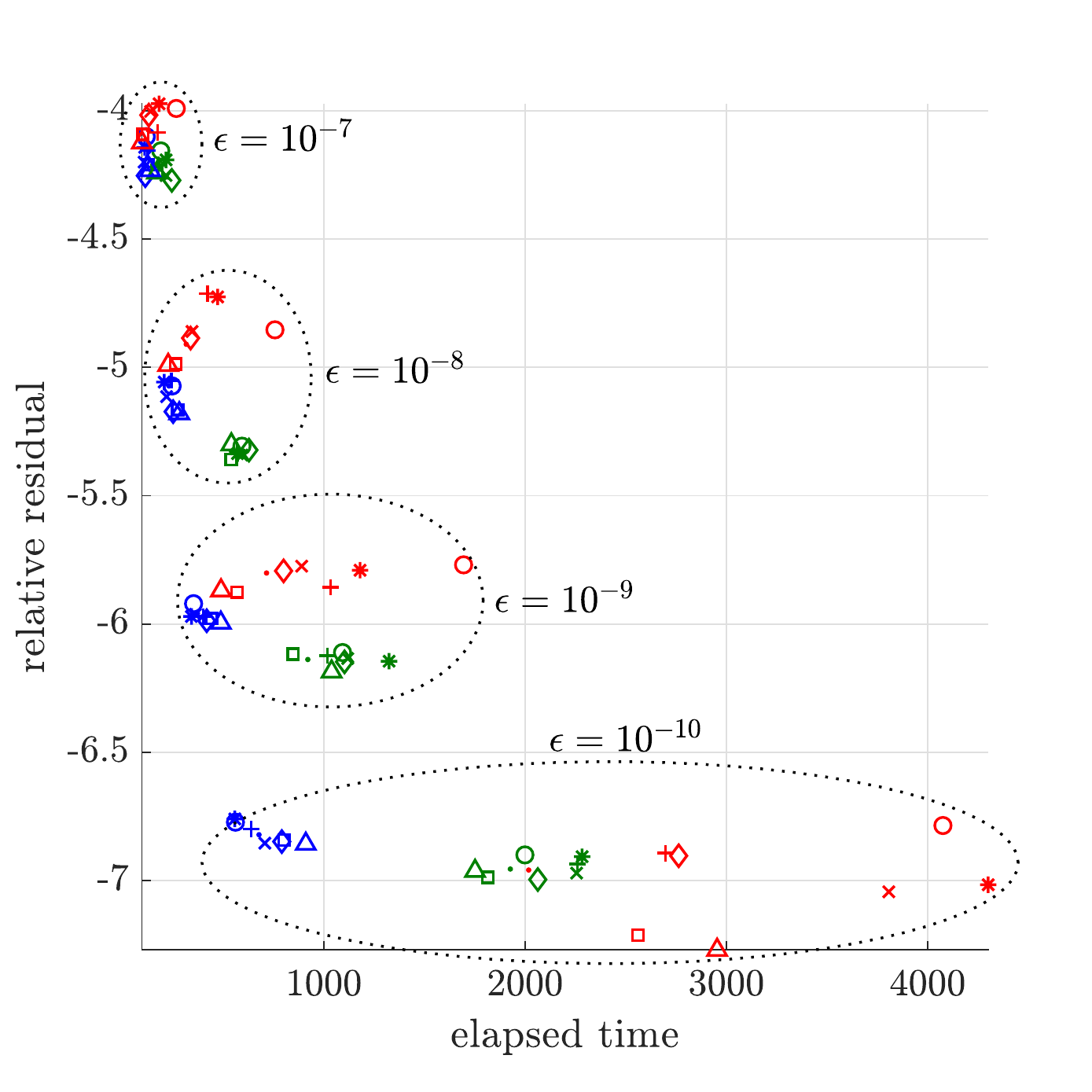}}
\subfloat[$m=24$, {[exp2]}]  {\label{fig:timing_hard2b}\includegraphics[scale=.45]{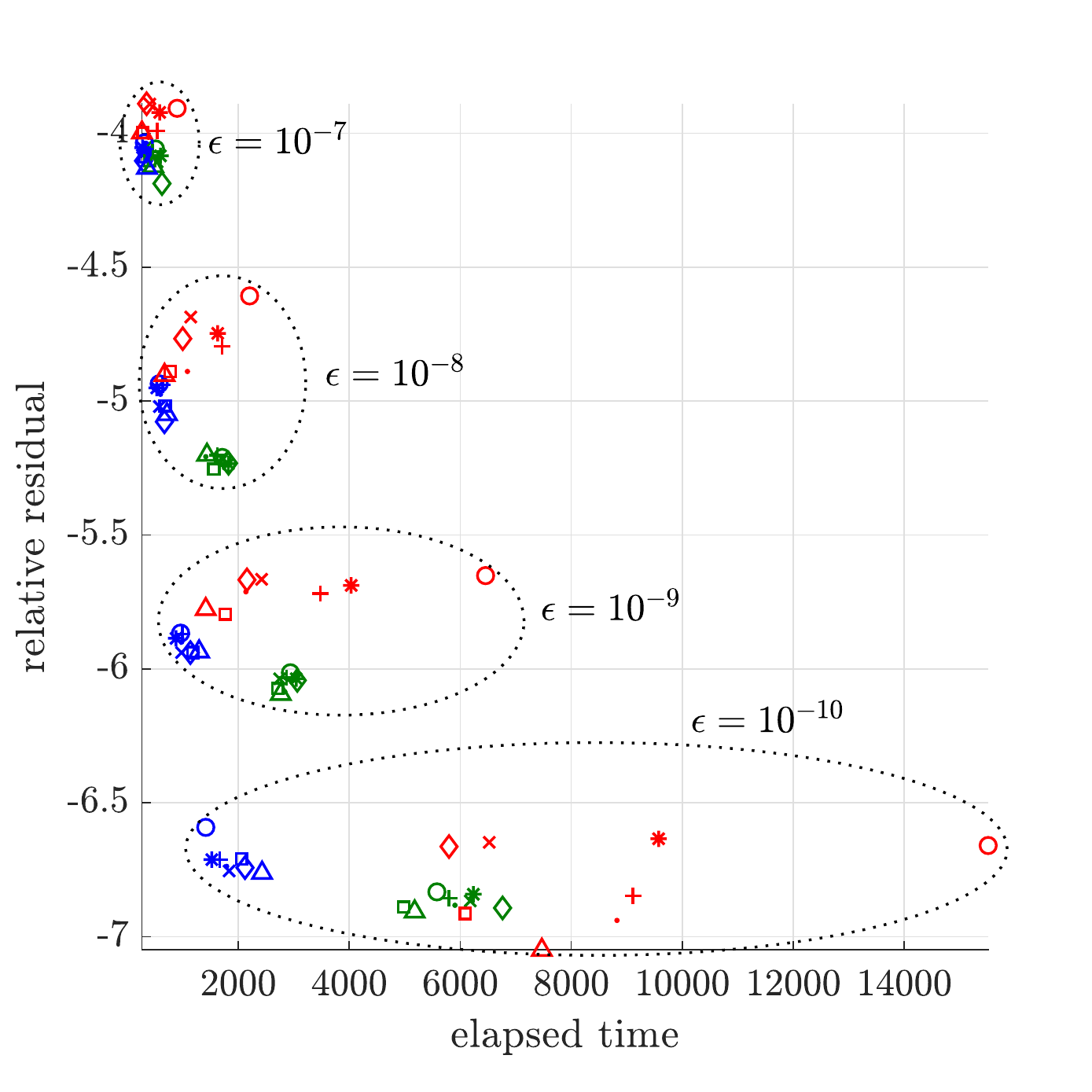}} \\
\caption{Computational timings (in seconds) of three EnhancedAEM methods for varying $k_{\max}$ and $n_{\text{update}}$. Timings of each method with each parameter set-up are averaged over 5 testing runs.}
\label{fig:timing_exp}
\end{figure}

In Figure \ref{fig:timing_exp}, we plot elapsed time (in seconds) against relative residual error for both
[exp1] and [exp2]. Note that the relative residual is computed afterwards in a post-processing step. Recall that the stopping condition for the outer iteration (see Algorithm \ref{alg:check_conv}) is not based on the relative residual (as this is expensive to compute). %{\color{red} We empirically found that upon termination, the relative residual error is around three orders of magnitude higher than the value of $\epsilon$ used}. 
The values of $\epsilon$ used for the stopping test for these results (see Algorithm \ref{alg:check_conv}) are
shown in the figure.  Note that for these experiments, the relative residual error is approximately
three orders of magnitude larger than $\epsilon$.
Results obtained with the EnhancedAEM methods with \pgdupdatename, \GSnameAbbrv, and \RstagenameAbbrv\ are
marked in red, green, and blue, respectively, and each configuration of $n_{\text{update}}$ and $k_{\max}$ is marked
with a different symbol. It can be seen from the figures that  
\begin{itemize}
\item the costs of \RstagenameAbbrv\ and \GSnameAbbrv\ are less sensitive to $n_{\text{update}}$ and $k_{\max}$ than those of \pgdupdatename;
\item \RstagenameAbbrv\ is more efficient for smaller values of $n_{\text{update}}$ whereas \GSnameAbbrv\ and
\pgdupdatename\ are better with larger $n_{\text{update}}$;
\item for \pgdupdatename\ and \GSnameAbbrv, relatively large
$n_{\text{update}} >10$ and $k_{\max}=2$ results in better performances, and,
for \RstagenameAbbrv, relatively small $n_{\text{update}} \leq 10$ and $k_{\max}=1$ results
in better performances.
\end{itemize}

Table \ref{tab:pvalues} reports the number of outer iterations $p$ required to achieve the stopping tolerance $\epsilon$ for problems [exp1] and [exp2] when \pgdupdatename, \GSnameAbbrv, and $\text{R-stage-}p$ are used. The benefit of using \RstagenameAbbrv\ becomes more pronounced as we seek highly accurate solutions with smaller $\epsilon$. 
Our general observation is that among the three enhancement approaches, the R-stage-p method is less sensitive
to choice of algorithm parameter inputs, scales better for larger problem sizes, and is the most effective of the three
approaches.

\begin{table}[tbhp]
{\footnotesize
  \caption{The number of outer iterations $p$ required to achieve the stopping tolerance $\epsilon$ for solving the problems [exp1] and [exp2] when \pgdupdatename, \GSnameAbbrv, and \RstagenameAbbrv\ are used. The reported values of $p$ are computed by averaging values of $p$ obtained with the eight different combinations of $n_\text{update}$ and $k_{\max}$ shown in the legend of Figure \ref{fig:timing_exp}.}\label{tab:pvalues}
\begin{center}
\renewcommand{\arraystretch}{1.25}
  \begin{tabular}{|c|c|c|c||c|c|c|} 
\multicolumn{7}{c}{[exp1]}\\
\hline
& \multicolumn{3}{c||}{$m=20$} &\multicolumn{3}{c|}{$m=24$}\\
\hline 
& \pgdupdatename\ & \GSnameAbbrv\ & \RstagenameAbbrv\ & \pgdupdatename\ & \GSnameAbbrv\ & \RstagenameAbbrv\ \\
\hline
$\epsilon = 10^{-7}$ &163.8 & 160.4 & 152.9 & 184.9 & 177.8 & 173.0\\
$\epsilon = 10^{-8}$ & 264.6 &	273.9 & 259.5 & 306.6 & 312.3 & 296.7\\
$\epsilon = 10^{-9}$ & 356.3 & 363.7 & 340.1 & 415.0 & 421.5 & 397.3\\
$\epsilon = 10^{-10}$ & 531.1 & 520.6 & 486.0 & 609.4 & 593.7 & 563.9\\
\hline
\multicolumn{7}{c}{}\\
\multicolumn{7}{c}{[exp2]}\\
\hline
& \multicolumn{3}{c||}{$m=20$} &\multicolumn{3}{c|}{$m=24$}\\
\hline 
& \pgdupdatename\ & \GSnameAbbrv\ & \RstagenameAbbrv\ & \pgdupdatename\ & \GSnameAbbrv\ & \RstagenameAbbrv\ \\
\hline
$\epsilon = 10^{-7}$ &293.1&287.7&282.1&344.0&334.9&330.6\\
$\epsilon = 10^{-8}$ &414.6&422.7&397.7&492.8&506.7&478.3\\
$\epsilon = 10^{-9}$ & 569.8 & 544.6 & 511.6 & 673.7 & 640.5 & 616.7\\
$\epsilon = 10^{-10}$ & 821.6 & 716.4 & 677.1 & 933.1 & 848.1 & 810.1\\
\hline
\end{tabular}
\end{center}
}
\end{table}

We now briefly consider a second benchmark problem whose solution matrix has different rank characteristics and for which low-rank solvers ought to perform well.

\subsection{Benchmark problem 2: fast decay coefficients}
We define the random field $a(x,\xi)$ as in \eqref{affine} but now we choose $\xi_{i} \sim U(-1,1)$ and the functions $a_{i}(x)$ have coefficients that decay more rapidly than in the first benchmark problem. The details of this problem can be found in \cite{eigel2017adaptive}.
Specifically, the coefficients of the expansion are
\begin{equation*}
a_{0}=1, \qquad a_i(x) = \alpha_i \cos (2\pi \varrho_1 (i) x_1) \cos (2\pi \varrho_2(i) x_2), \quad i=1,2,\ldots, m
\end{equation*}
where $\alpha_i =\bar{\alpha}i^{-\sigma}$ with $\sigma>1$ and $\bar{\alpha}$ satisfies  $0<\bar{\alpha}<1/\zeta(\sigma)$, where $\zeta$ is the Riemann zeta function. Furthermore, $\varrho_1(i) = i - k(i) (k(i)+1)/2$ and $\varrho_2(i) = k(i) - \varrho_1(i)$ where $k(i)=\lfloor -1/2 + \sqrt{1/4+2i} \rfloor$. Our implementation is based on the \textsc{Matlab} software package \textsc{S-IFISS} \cite{Silvester2016}. In the following experiment, we choose $\sigma=4$ and $\bar{\alpha}=0.832$. The parameter $\sigma$ controls the rate of algebraic decay of the coefficients. The specific choice $\sigma=4$ leads to fast decay and this causes the true solution matrix to have a lower rank than in the first benchmark problem.

\begin{figure} %[!h]
\begin{minipage}{.45\linewidth}
\centering
\includegraphics[scale=.5]{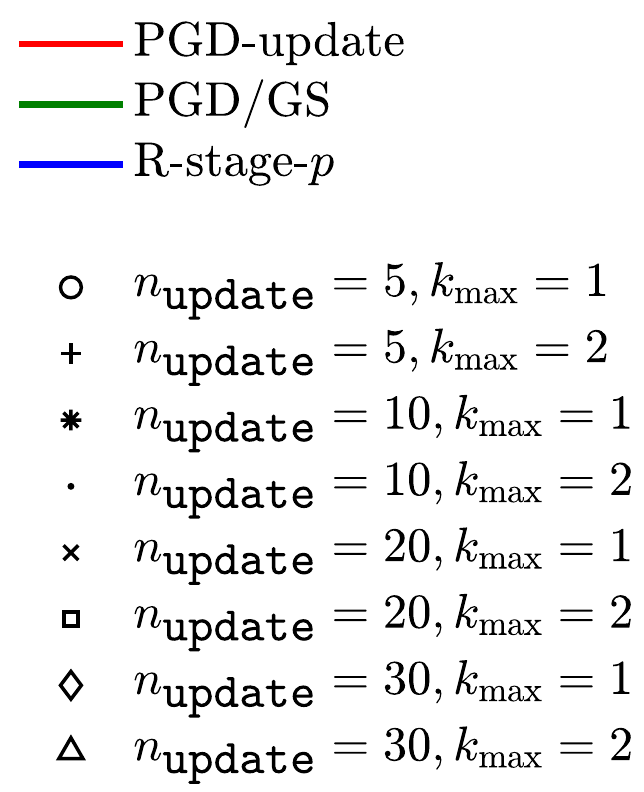}
\end{minipage}
\begin{minipage}{.5\linewidth}
\centering
%\subfloat[$m=20$]  {\label{fig:timing_fast}\includegraphics[scale=.45]{figs/[FAST]M20_scatter_all}}
\includegraphics[scale=.45]{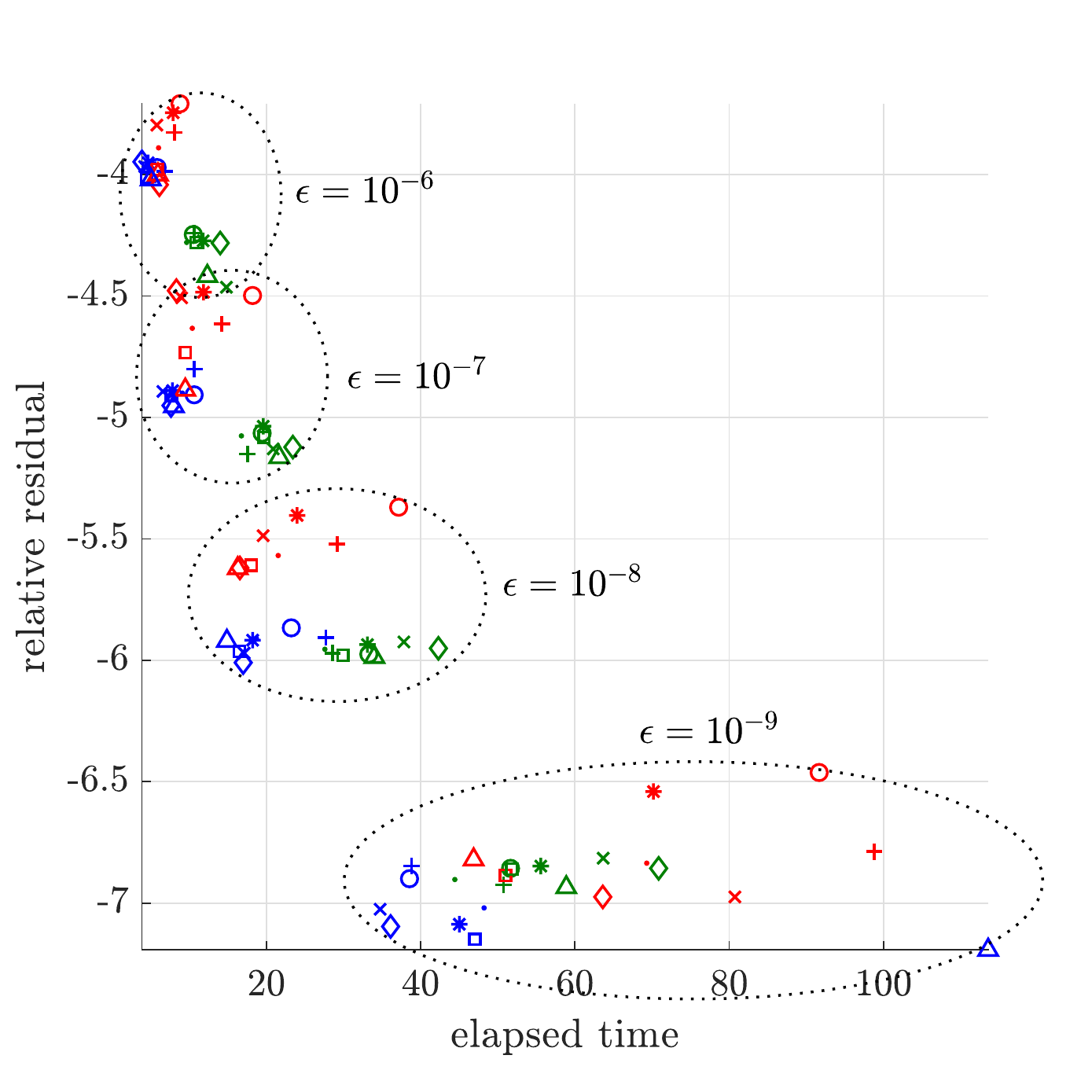} 
\end{minipage}
\caption{Computational timings (in seconds) of three EnhancedAEM methods for varying $k_{\max}$ and $n_{\text{update}}$. Timings of each method with each parameter set-up are averaged over 5 testing runs.}
\label{fig:timing_exp_fast}
\end{figure}

We investigate computational timings of the EnhancedAEM methods with the same experimental settings used in Section \ref{sec:time}. Here, we vary the stopping tolerance for the outer iterations as $\epsilon =\{10^{-9},10^{-8},10^{-7}\,10^{-6}\}$ and we choose the same values of $n_\text{update}$ and $k_{\max}$ as before.  Figure \ref{fig:timing_exp_fast} reports elapsed time (in seconds) against relative residual error.  In nearly all cases, our observations agree with the findings in Figure \ref{fig:timing_exp}. However, the impact of $n_\text{update}$ is slightly less clear for these tests. The R-stage-p method is generally still less sensitive than the other two methods to the
choices of $n_\text{update}$ and $k_{\max}$, with one exception, indicated by the blue triangle marker, which is located to the far right in Figure \ref{fig:timing_exp_fast}. With $n_\text{update}=30$, $k_{\max}=2$, and $\epsilon=10^{-9}$ (giving the right-most blue triangle), the \RstagenameAbbrv\ method does not meet the stopping criterion until $p\approx 125$, which is larger
than the value $p\approx 90$ needed for the other choices of algorithm inputs. We attribute
this to the large number of steps (30) between enhancements; in this case, the method fell just short of the stopping criterion after 90 steps.
Finally, we report the number of outer iterations $p$ required to achieve the stopping tolerance $\epsilon$ in Table \ref{tab:pvalues_fast}. As the true solution matrix has an intrinsic low-rank structure, the reported values of $p$ are much smaller than those shown in Table \ref{tab:pvalues}.
\begin{table}[!t]
{\footnotesize
  \caption{The number of outer iterations $p$ required to achieve the stopping tolerance $\epsilon$ for solving the second benchmark problem when \pgdupdatename, \GSnameAbbrv, and \RstagenameAbbrv\ are used. The reported values of $p$ are computed by averaging values of $p$ obtained with the eight different combinations of $n_\text{update}$ and $k_{\max}$ shown in the legend of Figure \ref{fig:timing_exp_fast}.}\label{tab:pvalues_fast}
\begin{center}
\renewcommand{\arraystretch}{1.25}
  \begin{tabular}{|c|c|c|c|} 
\hline 
& \pgdupdatename\ & \GSnameAbbrv\ & \RstagenameAbbrv\ \\
\hline
$\epsilon = 10^{-6}$ &43.7 & 49.0 & 30.1 \\
$\epsilon = 10^{-7}$ & 58.3 & 68.3 & 41.4 \\
$\epsilon = 10^{-8}$ & 81.7 & 91.7&	61.3 \\
$\epsilon = 10^{-9}$ & 130.9 & 121.6 & 91.6 \\
\hline
\end{tabular}
\end{center}
}
\end{table}

\subsection{Further Extensions} We also tested all the AEM methods on matrix equations obtained from stochastic Galerkin finite element discretizations of stochastic convection-diffusion problems \cite[Section 5.2]{lee2017preconditioned}, where the randomness is in the diffusion coefficient as in Section \ref{sec:ex_spd}. Although the energy norm cannot be defined for this problem because it has a non-symmetric operator, the same projection framework described herein can be applied to compute approximate solutions. Experiments (not reported here) were conducted similar to the ones in Sections \ref{sec:svs}--\ref{sec:acc}. We observed that the proposed \RstagenameAbbrv\ method produces qualitatively better approximate factors $V_{p}$ and $W_{p}$, as measured in the error metrics used in Sections \ref{sec:svs}--\ref{sec:acc}, than the \baselinename\ AEM method and the other two EnhancedAEM methods. 

\section{Conclusions} \label{sec:conclusions}

In this study, we have investigated several variants of alternating minimization methods to compute
low-rank solutions of linear systems that arise from stochastic Galerkin finite element
discretizations of parameterized elliptic PDEs.  Using a general formulation of alternating
energy minimization methods derived from the well-known general projection method, our starting
point was a variant of the stagewise ALS method, a technique for building rank-$p$ approximate
solutions developed for matrix completion and matrix sensing.  Our main contribution consists of a
combination of this approach with so-called enhancement procedures of the type used for PGD methods \cite{nouy2007generalized, nouy2010proper}
in which rank-one approximate solutions are enhanced by adaptive use of higher-rank
quantities that improve solution quality but limit costs by adaptively restricting the rank of
updates.  
Experimental results demonstrate that the proposed \GSnameAbbrv\ and \RstagenameAbbrv\ methods
produce accurate low-rank approximate solutions built from good approximations
of the singular vectors of the matricized parameter-dependent solutions.
Moreover, the results show that the $\text{R-stage-}p$ method scales better for larger
problems, is less sensitive to algorithm inputs, and produces approximate solutions
in the fastest times.

\section{Acknowledgements} \label{sec:ack}
This paper describes objective technical results and analysis. Any subjective
views or opinions that might be expressed in the paper do not necessarily
represent the views of the U.S. Department of Energy or the United States
Government. Sandia National Laboratories is a multimission laboratory managed
and operated by National Technology and Engineering Solutions of Sandia, a
wholly owned subsidiary of Honeywell International, for the U.S. Department of
Energy's National Nuclear Security Administration under contract DE-NA-0003525. This work was supported by the U.S. Department of Energy Office of Advanced Scientific Computing Research, Applied Mathematics program under award DE-SC0009301 and by the U.S. National Science Foundation under grant DMS1819115.

\bibliographystyle{siamplain}
\bibliography{ref}
\end{document}

%% file: ex_shared.tex
\usepackage{lipsum}
\usepackage{amsfonts}
\usepackage{graphicx}
\usepackage{epstopdf}
\usepackage{algpseudocode}
\ifpdf
  \DeclareGraphicsExtensions{.eps,.pdf,.png,.jpg}
\else
  \DeclareGraphicsExtensions{.eps}
\fi

\newsiamremark{remark}{Remark}
\newsiamremark{hypothesis}{Hypothesis}
\crefname{hypothesis}{Hypothesis}{Hypotheses}
\newsiamthm{claim}{Claim}

\headers{AEM Methods for Matrix Equations}{K. Lee, H. C. Elman, C. E. Powell, and D. Lee}

\title{Alternating Energy Minimization Methods for Multi-term Matrix Equations}

\author{Kookjin Lee\thanks{Extreme Scale Data Science and Analytics Department, Sandia National Laboratories %, Livermore 
(\email{koolee@sandia.gov}).}
\and Howard C. Elman\thanks{Department of Computer Science, University of Maryland, College Park (\email{elman@cs.umd.edu}).}
\and Catherine E. Powell\thanks{Department of Mathematics, University of Manchester,  UK (\email{c.powell@manchester.ac.uk}).}
\and Dongeun Lee\thanks{Department of Computer Science, Texas A\&M University-Commerce %, %Commerce, TX 75428 
(\email{dongeun.lee@tamuc.edu}).}}

\usepackage{amsopn}